\documentclass{article}
\usepackage{graphicx}
\usepackage{amsfonts}
\usepackage{hyperref}

\newcommand{\goesto}{\mbox{$\rightarrow$}}

\newcommand{\la} {\langle}
\newcommand{\ra} {\rangle}

\newcommand{\QED} {\rule{1mm}{2.5mm}}
\newtheorem{definition}{Definition}
\newtheorem{Theorem}{Theorem}
\newtheorem{theorem}{Theorem}

\newtheorem{lemma}[Theorem]{Lemma}

\newtheorem{corollary}[Theorem]{Corollary}

\newcommand{\Euc}{\mathbb{E}}
\newcommand{\Real}{\mathbb{R}}

\newcommand{\Topo} {\mbox{$\mathbb{T}$}}
\newcommand{\Reg}{\mbox{$\mathcal R$}}
\newcommand{\Convex}{\mbox{$\mathcal C$}}

\title{Metric topologies over some categories of simple open regions
in Euclidean space}
\author{Ernest Davis
\\ Dept. of Computer Science \\ New York University
\\ davise@cs.nyu.edu}

\begin{document}
\maketitle
\begin{abstract}
What does it mean for a shape to change continuously? Over the space of
convex regions, there is only one ``reasonable'' answer. However, over a 
broader class of regions, such as the class of star-shaped regions, there
can be many different ``reasonable'' definitions of continuous shape change.

We consider the relation between topologies induced by a number of metrics
over a number of limited categories of open bounded regions in $\Euc^{n}$.
Specifically, we consider a homeomorphism-based metric; the 
Hausdorff metric; the dual-Hausdorff metric; the symmetric difference
metric; and the family of Wasserstein metrics; and the topologies that they
induce over the space of convex regions; the space of convex regions and unions of two
separated convex regions; and the space of star-shaped regions. We
demonstrate that:
\begin{itemize}
\item Over the space of convex regions, all five metrics, and indeed any
metric that satisfies two general well-behavedness constraints, induce
the same topology.
\item Over the space of convex regions and unions of two
separated convex regions, these five metrics
are all ordered by ``strictly finer than" relations. In descending order
of fineness, these are: 
the homeomorphism-based, the dual-Hausdorff, the Hausdorff,
the Wasserstein, and the symmetric difference. Also, Wasserstein metrics
are strictly ordered among themselves.
\item  Over the space of star-shaped regions, the topologies induced
by the Hausdorff metric, the
symmetric-difference metric, and the Wasserstein metrics are incomparable
in terms of fineness.
\end{itemize}
\end{abstract}


\section{Introduction}
\label{secIntro}
In many applications in physical reasoning and in computer graphics, shapes
deform continuously. However,  what kinds of functions from time to shapes
count as ``continuous'' depends on the topology of the space
of regions; and this, as we will discuss here, is not as clear-cut as one
might suppose.

Over the space of points in $\Euc^{n}$, there are a number of different metrics
in common use: the standard Euclidean distance, the Manhattan distance, and
more generally the Minkowski distance with parameter $p$. 
But all of these, except 
the discrete metric, are fundamentally similar, in the sense that they generate
the same topology.  If a sequence ${\bf x}_{1}, {\bf x}_{2}$ 
converges to $\bf y$
in any of them, then it converges to $\bf y$ in all of them; and if a function 
$\phi({\bf x})$ is continuous in any of them, it is continuous in all of them.

When one considers the space of regions in $\Euc^{n}$, however, the situation
is very different. Here, again, there are many different possible 
natural metrics,
with no obvious clear favorite, and these are fundamentally different in the
sense that they generate different topologies 
(Davis 2001; Galton 2000).

\begin{figure}
\begin{center}
\includegraphics[width=4in]{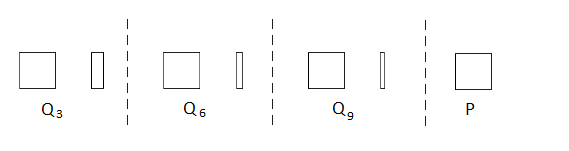}
\caption{A sequence of regions that converges in the symmetric-difference
metric but not in the Hausdorff metric}
\end{center}
\label{figAreaVsHausdorff}
\end{figure}

For instance, 
(figure~\ref{figAreaVsHausdorff}) 
consider the sequence of regions in the plane
${\bf Q}_{1}, {\bf Q}_{2} \ldots$ where 
${\bf Q}_{i}=((0,1) \times (0,1)) \cup ((2,2+1/i) \times (0,1))$. Let
${\bf P} = (0,1) \times (0,1)$. 
If one measures the difference
between two regions {\bf X} and {\bf Y} as the area of their symmetric difference 
\[V({\bf X,Y}) = \mbox{area}(({\bf X} \setminus {\bf Y}) \: \cup \:
({\bf Y} \setminus {\bf X})) \]
then $V({\bf Q}_{i},{\bf P}) = 1/i$, so the sequence 
${\bf Q}_{1}, {\bf Q}_{2}$ converges to ${\bf P}$. If one measures it
using the Hausdorff distance $H({\bf X},{\bf Y})$, then
$H({\bf Q}_{i},{\bf P}) > 1$ for all $i$,
so the sequence does not converge to {\bf P}.

In this paper, we consider limited classes of regions and well-known 
metrics that satisfy two specified well-behavedness conditions. We consider the
relations between the topologies that these metrics generate over these classes.
We prove results of two general flavors. 
First, in section~\ref{secConvex}, we show that,
over the space of {\em convex\/} regions 
there is only one natural metric topology. More precisely, the theorem 
show that any metric
satisfying these well-behavedness conditions generate the same topology. 
Thus, for
instance, there is no way to construct an example analogous to 
figure~\ref{figAreaVsHausdorff} using convex regions; if a sequence of convex
open regions converges to a convex open region
in the area metric, it also converges in the Hausdorff
metric, and in any other well-behaved metric over regions.

The second flavor of result, show that, as
figure~\ref{figAreaVsHausdorff} illustrates, if one expands the space of regions
under consideration to a broader class of
regions,  then the different metrics we consider generate 
different topologies.

Section~\ref{secNotation} will introduce notational conventions and basic
functions.
Section~\ref{secWellBehaved} will define our two well-behavedness conditions:
A well-behaved topology ``supports continuous morphing''and 
``satisfies the region
separation condition''.
Section~\ref{secMetrics} defines the metrics we will consider: \\
\hspace*{2em} A homeomorphism-based metric $M({\bf A,B})$; \\
\hspace*{2em} The Hausdorff metric $H({\bf A,B})$; \\
\hspace*{2em} The dual-Hausdorff metric $H^{d}({\bf A,B})$; \\
\hspace*{2em} The symmetric-difference metric $V({\bf A,B})$; and \\
\hspace*{2em} The family of Wasserstein metrics $W^{\psi}{\bf A,B})$

We demonstrate that:
\begin{itemize}
\item Over the space of convex regions, all five metrics, and indeed any
metric that satisfies two general well-behavedness constraints, induce
the same topology (section~\ref{secConvex}).
\item Over the space of convex regions and unions of two
separated convex regions, these five metrics
are all ordered by ``strictly finer than" relations. In descending order
of fineness, these are: 
the homeomorphism-based, the dual-Hausdorff, the Hausdorff,
the Wasserstein, and the symmetric difference. Also, Wasserstein metrics
are strictly ordered among themselves (section~\ref{secTwoConvex}).
\item  Over the space of star-shaped regions, the topologies induced
by the Hausdorff metric, the
symmetric difference metric, and the Wasserstein metrics are incomparable 
in terms of fineness (section~\ref{secStar}).
\end{itemize}

\subsection{Notation and basic concepts} 
\label{secNotation}

$\Real$ is the space of real numbers.
$\Euc^{n}$ is $n$-dimensional Euclidean
space. We will generally assume that $n \geq 2$; many of our concepts become
vacuous or trivial in one-dimensional space, though some carry over. 

Real numbers and distances will be notated with italicized variables: $x,d$.

Points in $\Euc$ will be notated with boldface lower-case variables: ${\bf p,q}$.
In some of the proofs, it will be convenient to choose an origin and notate
points as vectors: $\vec{p}, \vec{q}$.
The standard Euclidean 
distance between points {\bf x} and {\bf y} will be denoted $d({\bf p,q})$.

Subsets of $\Euc^{n}$ will be notated with boldface capital letters: $\bf P,Q$.

A {\em region}
will be a open subset of $\Euc^{n}$ that is bounded and equal to the 
interior of its closure (topologically regular).
The class of all regions in $\Euc^{n}$ will be denoted $\Reg$ (the
dimension of the space being left implicit).
The closure of {\bf A} is denoted $\bar{\bf A}$.
The topological boundary of region {\bf A} (i.e. the closure of 
{\bf A} minus {\bf A}) is denoted  @{\bf A} = $\bar{\bf A} \setminus {\bf A}$.

The $n$-dimensional volume of region {\bf A} is denoted $v({\bf A})$.

The open ball of radius $d$ centered
at point ${\bf p}$ is denoted ${\bf B}({\bf p},d) \subset \Euc^{n}$. 

The {\em radius} of region {\bf A} at point ${\bf o}  \in {\bf A}$
is the radius of the largest spherical 
open ball that fits inside ${\bf A}$. 
The radius (1 argument) of region {\bf A} is its maximal radius.
$\mbox{radius}({\bf A})=\max_{{\bf o} \in {\bf A}} \: \mbox{radius}({\bf A,o})$.

The {\em diameter} of {\bf A} is the maximal distance between
two points in $\bar{\bf A}$:
$\mbox{diameter}({\bf A}) = \sup_{{\bf p,q} \in \bar{A}} d({\bf p,q})$.

The {\em distance from point ${\bf p}$ to region {\bf Q}} is the distance
from {\bf p} to the closest point in the closure of {\bf Q}.
\[ d({\bf p},{\bf Q}) = \min_{{\bf q} \in \bar{\bf Q}} \: d({\bf p,q}) \]

The {\em distance} between regions {\bf A} and
{\bf B} is the smallest distance between points in their closure:
$d({\bf A,B}) = \min_{{\bf a}\in \bar{\bf A},{\bf b} \in \bar{\bf B}}
d({\bf a,b})$. The distance $d({\bf A,B})$ is not, of course, a metric over
regions.

\begin{definition}
Let {\bf P} be a region. Let $\delta > 0$. 

The {\em dilation} of {\bf P}
by $\delta$ is the set of all points within $\delta$ of {\bf P}. \\
dilate(${\bf P},\delta$) =
$\{ {\bf w} \: | \: d({\bf w,P}) \leq \delta \}$. 

The {\em erosion\/} of {\bf P} by $\delta$ is the set of all points more than 
$\delta$ from the complement of {\bf P}. \\
erode(${\bf P},\delta$) =
$\{ {\bf x} \: | \: d({\bf x,P}^{c}) \geq \delta \}$. 

The {\em outer shell} of {\bf P} by $\delta$, ${\bf O}({\bf P},\delta) =
\mbox{dilate}({\bf P},\delta) \setminus {\bf P}$.

The {\em inner shell} of {\bf P} by $\delta$, ${\bf I}({\bf P},\delta) =
{\bf P} \setminus \mbox{erode}({\bf P},\delta)$.
\end{definition}

The regularization of $\bf X \subset \Euc^{n}$ is the interior of the closure 
of $\bf X$. Boolean operators, as applied to regions, are implicitly regularized.
For instance if 
${\bf P}=(0,1) \times (0,1)$,
${\bf Q}=(1,2) \times (0,1)$, and
${\bf R}=(0,2) \times (0,1)$, then
${\bf P}  \cup {\bf Q} = {\bf R}$ and ${\bf R} \setminus {\bf Q}={\bf P}$.

Subsets of $\Reg$ -- that is, sets  of subsets of $\Euc^{n}$ --- 
will be denoted
using calligraphic letters: $\mathcal U$, $\mathcal V$.  \\
In particular $\mathcal C$ is the set of all convex regions. \\
${\mathcal D}^{2}$ is the set of all regions that are the union of two 
separated convex regions: \\
${\mathcal D}^{2} = \{ {\bf X} \cup {\bf Y} \: | \:
{\bf X,Y} \in {\mathcal C}, d({\bf X,Y}) > 0)$. \\
$\mathcal D$ will be the set of all regions that are either a single convex
region or the union of two separated convex regions; thus
${\mathcal D} = {\mathcal C} \cup {\mathcal D}^{2}$. \\
${\mathcal S}$ will be the set of all bounded, star-shaped regions.

We will use 
$\mu: \Reg \times \Reg \mapsto \Real$ 
to represent a generic metric 
over $\Reg$; that is $\mu({\bf A,B})$ is some measure of the
difference between regions $\bf A$ and $\bf B$  that satisfies the
standard axioms for metrics.
We will use upper-case italic letters for specific
metrics, as defined in section~\ref{secMetrics}; for instance, the Hausdorff
distance is denoted $H({\bf P},{\bf Q})$.

Otherwise, the font of function symbols will correspond to the type of the
value returned by the function. In particular, 
the ball of radius $d$ relative to the metric $\mu$ centered at region $\bf P$
is denoted 
${\mathcal B}_{\mu}({\bf P},d) = \{ {\bf Q} \: | \: \mu({\bf P},{\bf Q}) < d \}$

Finally $\Topo_{\mu}$ will be the topology generated by metric $\mu$
over $\Reg$; since a topology is a set of open sets, 
$\Topo_{\mu}$ is a set of sets of subsets of $\Euc^{n}$.

Throughout this paper, the phrases 
``$\Topo_{\alpha}$ is finer than $\Topo_{\beta}$''
or ``is coarser'', if unqualified, are to be interpreted as a non-strict relation;
that is, as ``finer than or equal to'' or ``coarser than or equal to''. When a
strict relation is intended, the phrases ``strictly finer/coarser'' will be used.
The phrase ``$\Topo_{\alpha}$ is not finer/coarser than $\Topo_{\beta}$'' will
mean ``It is not the case that $\Topo_{\alpha}$ is finer/coarser than 
$\Topo_{\beta}$.''

\section{Well-behaved topologies}
\label{secWellBehaved}
\begin{definition}
Let $\mathcal U$ be a set of regions (a subset of $\Reg$).
A {\em history} over $\mathcal U$ is a function $\phi: [0,1] \mapsto \mathcal U$.
\end{definition}

\begin{definition}
A {\em morphing} over $\Euc^{n}$ is a  uniformly continuous
function $\psi : [0,1] \times \Euc^{n} \mapsto \Euc^{n}$ 
with the following properties:
\begin{itemize}
\item[a.] $\psi(0,\cdot)$ is the identity over $\Euc^{n}$
\item[b.] For $t \in [0,1]$, $\phi(t, \cdot$) is a homeomorphism of $\Euc^{n}$
to itself.
\end{itemize}
\end{definition}

\begin{definition}
\label{defCorrespondsMorphing}
A history $\phi:\Real \mapsto \Reg$ {\em corresponds to morphing $\psi$} if
$\phi(t) = \psi(t,\phi(0))$.
\end{definition}

\begin{definition}
A topology $\Topo$ over a subspace $\mathcal U$ of $\Reg$ {\em supports 
continuous morphing} if every history over $\mathcal U$ that 
corresponds to a morphing is continuous  relative to $\Topo$.
\end{definition}

Intuitively, if you start with a spatial region {\bf A} and you morph it around
continuously relative to the regular spatial topology, then its trajectory
as a function of time is continuous in $\Topo$. This is an upper bound on
the fineness of $\Topo$; the topology cannot be so fine that morphings are
discontinuous. If $\Topo$ supports continuous morphing and $\Topo'$ is coarser
than $\Topo$, then $\Topo'$ also supports continuous morphing.

The following is an example of a metric that does not support 
continuous morphing.
Let $\mathcal U$ be the set of regions in $\Euc^{2}$ with a finite perimeter.
Define the metric over $\mathcal U$, 
$\mu({\bf X,Y}) = H({\bf X,Y}) + 
|\mbox{perimeter}({\bf X})-\mbox{perimeter}({\bf Y})|$.
Then one can easily define a morphing in which $\phi(0)$ is the unit square
and $\phi(t)$ is the unit square with a saw-toothed boundary, where the
teeth are at $45^{\circ}$ and the length of the teeth is $t$. Then
for all $t > 0$, the perimeter of $\phi(t)$ is approximately $4 \sqrt{2}$,
so the morphing is not continuous relative to $\Topo_{\mu}$.

\begin{definition}
\label{defSeparates}
A topology $\Topo$ over $\Reg$ {\em satisfies the region separation 
condition\/}
if the following hold for any regions ${\bf P,Z} \in \Reg$:
\begin{itemize}
\item[i.] If ${\bf P} \cap {\bf Z} = \emptyset$, 
then in $\Topo$ there exists a neighborhood $\mathcal U$ of $\bf P$
such that no superset of ${\bf Z}$ is in $\mathcal U$.

\item[ii.] If ${\bf P} \supset {\bf Z}$, 
then in $\Topo$ there exists a neighborhood $\mathcal U$ of $\bf P$
such that no region that is disjoint from ${\bf Z}$ is in $\mathcal U$.
\end{itemize}
\end{definition}

\begin{lemma}
\label{lemSuddenEmergence}
Let $\Topo$ be a topology over $\Reg$ that satisfies the region separation
condition.
Let $\phi:\Real \mapsto \Reg$ be a history that is
continuous under $\Topo$. Let
${\bf Z} \in \Reg$ be any open region. Then
there exists a neighborhood $U$ of 0 such that 
\begin{itemize}
\item if ${\bf Z} \cap \phi(0) = \emptyset$ then
there is no $t \in U$ such that ${\bf Z} \subset \phi(t)$; 
\item if ${\bf Z} \subset \phi(0)$ then
there is no $t \in U$ such that ${\bf Z} \cap \phi(t) = \emptyset$.
\end{itemize}
\end{lemma}

{\bf Proof:} Taking ${\bf P}=\phi(0)$, 
construct the set $\mathcal U$ to satisfy the conclusion of 
definition~\ref{defSeparates}.
Take $U = \phi^{-1}({\mathcal U})$. By continuity, $U$ is open and
by construction it satisfies the conditions of the theorem.

\begin{definition}
A topology is {\em well-behaved} if it supports continuous morphing
and satisfies the region separation condition.
\end{definition}

It is immediate from the definitions  that if a topology supports continuous
morphing, then every coarser topology does; and that if a topology satisfies the
region separation condition, then every finer topology does.

\section{Metrics on the space of regions}
\label{secMetrics}
In this paper, we primarily consider five metrics, or families of metrics, over the
space of regions:
a homeomorphism-based metric $M({\bf A,B})$;
the Hausdorff metric $H({\bf A,B})$;
the dual-Hausdorff metric $H^{d}({\bf A,B})$;
the symmetric-difference metric $V({\bf A,B})$;
and the family of Wasserstein metrics $W^{\psi}({\bf A,B})$

Some other metrics will be discussed in passing at various points.

\subsection{Homeomorphism-based metric}
\label{secHomeoMetrics}
There are a number of different ways of defining the difference between
two regions {\bf A} and {\bf B} in terms of homeomorphisms between them or
between their boundaries. Perhaps the oldest and the best known is the
Fr\'{e}chet distance. In this paper we will use 
the {\em homeomorphism distance}
$M({\bf A,B})$, defined as follows:

Let {\bf A} and {\bf B} be two regions in $\Euc^{n}$.  
Let $\Gamma({\bf A},{\bf B})$ be the set of all homeomorphisms $\gamma$ 
of $\Euc^{n}$ to itself such that $\gamma({\bf A})={\bf B}$.
Define the metric 
\[ M({\bf A,B}) = \inf_{\gamma \in \Gamma} \sup_{{\bf x}\in \Euc^{n}}
d({\bf x},\gamma({\bf x})) \]

(If $\Gamma=\emptyset$ --- that is, there are no homeomorphisms of the space
that map {\bf A} to {\bf B} --- then $M({\bf A,B}) = \infty$.)

In other words: for any $\gamma$ that is an homeomorphism from $\Euc^{n}$ to
itself and that maps $\bf A$ to $\bf B$, we define a cost which is the maximum
distance from $\bf x$ to $\gamma({\bf x})$ for any $\bf x$ in $\Euc^{n}$
We then define the metric $M({\bf A,B})$ as the 
smallest cost attained by any such $\gamma$ (more precisely, the
infimum). 

\begin{theorem}
\label{thmMContinuous}
The topology $\Topo_{M}$ supports continuous morphings over $\Reg$.
\end{theorem}

{\bf Proof:} Immediate from the definition.

A converse of theorem~\ref{thmMContinuous} would be the claim that 
if a history is continuous relative to $\Topo_{M}$ then it corresponds to a 
morphing. I suspect that this is true, but have not been able to prove it.

\subsection{The Hausdorff and dual-Hausdorff metrics}

The {\em one-sided Hausdorff distance\/} from region {\bf P} to {\bf Q}
is the supremum over points {\bf p} in {\bf P} 
of the distance from {\bf p} to {\bf Q}.

\[ H^{1}({\bf P,Q}) = \sup_{{\bf p} \in {\bf P}} d({\bf p,Q}) \]

The {\em Hausdorff distance\/} between regions and {\bf P} to {\bf Q}
is the maximum of 
(the one-sided Hausdorff distance from {\bf P} to {\bf Q}) and
(the one-sided Hausdorff distance from {\bf Q} to {\bf P}) 

\[ H({\bf P,Q}) = \max(H^{1}({\bf P,Q}), H^{1}({\bf Q,P})) \]

The {\em dual-Hausdorff distancei\/}  (Davis 1995)
is the maximum of (the Hausdorff distance between {\bf P} and {\bf Q})
and (the Hausdorff distance between the complements of {\bf P} and {\bf Q}).

\[ H^{d}({\bf P,Q}) = \max(H({\bf P,Q}), H({\bf Q}^{c},{\bf P}^{c})) \]

This metric is not discussed in (Deza and Deza 2006) but the proof that it is 
a metric over the
space of regular regions is immediate.

It is immediate from the definitions that for all regions, 
$H({\bf P,Q}) \leq H^{d}({\bf P,Q}) \leq M({\bf P,Q})$ and
therefore $\Topo_{M}$ is finer than $\Topo_{H^{d}}$ which is finer than
$\Topo_{H}$.

\begin{theorem}
\label{thmHausdorffContMorph}
Topologies $\Topo_{H^{d}}$ and $\Topo_{H}$ support continuous morphing
over $\Reg$.
\end{theorem}

{\bf Proof:} Immediate from theorem~\ref{thmMContinuous} together
with the above.

\begin{theorem}
\label{thmHausdorffSeparation}
The Hausdorff distance has the region separation property over $\Reg$.
\end{theorem}

{\bf Proof:} 
i. Let {\bf P, Z} be regions such that ${\bf P}\cap {\bf Z}
=\emptyset$. Let ${\bf Y} \supset {\bf Z}$. 
Let {\bf z} be a point in {\bf Z}. 
Then $H({\bf Y,P}) \geq d({\bf z,P})$. So for
$\epsilon < d({\bf z,P})$, the open ball
${\mathcal B}_{H}({\bf P},\epsilon)$ excludes all {\bf Z} and any 
superset of {\bf Z}.

ii. Let {\bf P, Z} be regions such that ${\bf Z} \subset {\bf P}$.
Let {\bf Y} be a region such that ${\bf Z}$ and {\bf Y} are disjoint.
Let {\bf z} be a point in {\bf Z}. 
Then $H({\bf Y,P}) \geq  \mbox{radius}({\bf P,z})$. So for
$\epsilon < \mbox{radius}({\bf Z,z})$, the open ball
${\mathcal B}_{H}({\bf P},\epsilon)$ excludes all {\bf Y} and any
subset of {\bf Y}.

\begin{corollary}
\label{corSeparation}
The metrics $M({\bf P,Q})$ and $H^{d}({\bf P,Q})$ 
have the region separation property over $\Reg$.
\end{corollary}

{\bf Proof:} It is immediate that, if a topology has the property, then any
finer topology also has the property.

\subsection{The symmetric-difference metric}
\label{secSymDiff}
Define the function {\bf S(P,Q)}: $\Reg \times \Reg \mapsto \Reg$ as the 
symmetric difference of regions {\bf P} and {\bf Q}: \\
${\bf S(P,Q)} = ({\bf P} \setminus {\bf Q})
\cup ({\bf Q} \setminus {\bf P})$

The {\em symmetric-difference} metric is the $n$-dimensional measure
of the symmetric difference: \\
$V({\bf P,Q}) = v({\bf S(P,Q}))$

\begin{theorem}
\label{thmTopoDHFinerVolume}
Over the space $\Reg$, $\Topo_{H^{d}}$ is finer than $\Topo_{V}$.
\end{theorem}

{\bf Proof:} See (Davis 2001), corollary 8.2.

\begin{theorem}
\label{thmVolumeSupportsMorphing}
$\Topo_{V}$ supports continuous morphings over $\Reg$.
\end{theorem}

{\bf Proof:} Immediate from theorem~\ref{thmHausdorffContMorph} and 
lemma~\ref{thmTopoDHFinerVolume}. 

\begin{theorem}
\label{thmVolumeSeparation}
$\Topo_{V}$ has the region separation property over $\Reg$.
\end{theorem}

{\bf Proof:} \\
i. Let {\bf P, Z} be regions such that ${\bf P}\cap {\bf Z}
=\emptyset$. Let ${\bf Y} \supset {\bf Z}$. 
Then ${\bf Z} \subset S({\bf P,Y})$,
$V({\bf P,Y}) \geq v({\bf Z})$. So for 
$\epsilon < v({\bf Z})$, the open ball
${\mathcal B}({\bf P},\epsilon)$ excludes 
{\bf Z} and any superset of {\bf Z}.

ii. Let {\bf P,Z} be regions such that ${\bf Z} \subset {\bf P}$.
Let {\bf Y} be a region such that ${\bf Z}$ and {\bf Y} are disjoint.
Then again ${\bf Z} \subset S({\bf P,Y})$,
So for $\epsilon < v({\bf Z})$, the open ball
${\mathcal B}({\bf P},\epsilon)$ excludes all sets disjoint from {\bf Z}.

\subsection{Wasserstein metrics}
\label{secWasserstein}
The family of  Wasserstein distances $W^{\psi}({\bf P,Q})$ are generalizations
of the ``earth-movers'' metric often used in comparing probability 
distributions.

{\bf Definition} A function $\psi:\Real^{\geq 0} \mapsto \Real^{\geq 0}$ is
a {\it Mulholland function} if it is continuous and monotonically increasing;
$\psi(0) = 0$; $\lim_{x \goesto \infty} \psi(x) = \infty$; and $\psi$ satisifies 
the Mulholland (1949) inequality 
\[ \psi^{-1}(\sum_{i=1}^{n} \psi(x_{i}+y_{i})) \leq \psi^{-1}(\sum_{i=1}^{n} \psi(x_{i}))+\psi^{-1}(\sum_{i=1}^{n} \psi(y_{i})) \]
The Minkowski inequality is the special case where $\phi(x)=x^{p}$.

The Wasserstein distance corresponding to a Mulholland function $\psi$ is a
metric over probability distributions. (It is usually defined using the particular
function $\psi(x) = x^{p}$. However, since the only property of $x^{p}$ that is
used in proving that the Wasserstein distance is a metric is that it satisfies
the Mulholland inequality, one can generalize it to use any Mulholland function
(Clement and Desch, 2008).)

\begin{definition}
\label{defWassersteinOverDist}
Let $\psi$ be a Mulholland function. Let $\theta({\bf x})$ and $\zeta({\bf x})$
be probability densities over $\Euc^{n}$. 
Let
$\gamma$ be a function
from $\Euc^{n}$ to $\Euc^{n}$ such that, if random
variable $X$ has density $\theta({\bf x})$ then 
$\gamma(X)$ will have density $\zeta({\bf x})$.
Define the integral
\[ I(\gamma) = \int_{{\bf x} \in \Euc^{n}} \theta({\bf x}) 
\cdot \psi(d({\bf x}, \gamma({\bf x}))) \: d{\bf x} \]

Let $\Gamma(\theta,\zeta)$ be the set of all such $\gamma$.
Then the {\em Wasserstein distance 
between $\theta$ and $\zeta$ corresponding to $\psi$} is defined as follows: 
\[ W^{\psi}(\theta,\zeta) = 
\inf_{\gamma \in \Gamma(\theta,\zeta)} \psi^{-1}(I(\gamma)) \]
\end{definition}

We adapt the above definition to be a distance between regions ${\bf P}$ and
${\bf Q}$ by taking $\theta$ and $\zeta$ to be the uniform distributions 
over ${\bf P}$ and ${\bf Q}$. 

\begin{definition}
For any region {\bf P}, $U_{\bf P}$ represents the uniform distribution
over {\bf P}: \\
$U_{\bf P}({\bf x}) = 1/v({\bf P})$ for ${\bf x} \in {\bf P}$. \\
$U_{\bf P}({\bf x}) = 0$ for ${\bf x} \not \in {\bf P}$. \\
\end{definition}

\begin{definition}
\label{defWassersteinReg1}
Let {\bf P} and {\bf Q} be regions in $\Reg$. Let $\psi$ be a Mulholland
function. Define $W^{\psi}({\bf P,Q})$ to be $W^{\psi}(U_{\bf P},U_{\bf Q})$
\end{definition}

We can reformulate this definition as follows:

\begin{definition}
Let $\bf P$ and $\bf Q$ be regions.
Let $\gamma$ be a function from
{\bf P} to {\bf Q}. We say that $\gamma$ is {\em uniform} if, for all 
${\bf X} \subset {\bf P}$,
$v(\gamma({\bf X}))) = v({\bf X})\cdot v({\bf Q})/v({\bf P})$. That is,
$\gamma$ preserves relative measure. 

Define the following two functions of $\gamma$ and ${\bf P}$:

\[ I^{\psi}(\gamma,{\bf P}) = 
\frac{1}{v({\bf P})} \cdot 
\int_{{\bf x} \in P} \psi(d({\bf x},\gamma({\bf x}))) \:
\mbox{d}{\bf x}  \]

\[ C^{\psi}(\gamma,{\bf P}) = 
\psi^{-1}(I^{\psi}(\gamma,{\bf P})) \]

Let $\Gamma({\bf P,Q})$ 
be the set of all uniform functions $\gamma$ from {\bf P} to
{\bf Q}. Then 
$W^{\psi}({\bf P,Q}) = 
\inf_{\gamma \in \Gamma({\bf P,Q})} C^{\psi}({\bf P,Q})$.
\end{definition}

In the case of the identity function $\psi(x)=x$, this can be given an
intuitive motivation as follows: Suppose that you have dirt uniformly spread
over $\bf P$ and you want to move it so that it is uniformly spread
out over $\bf Q$. To move a small piece of dirt of mass $m$ from 
{\bf x} to {\bf y} will cost $m \cdot d({\bf x,y})$. Then 
if you follow $\gamma$ as a guide for how to move the dirt, 
the total cost will be $C^{\psi}(\gamma)$. Thus
the cost of the cheapest way of moving the dirt is $W^{\psi}({\bf A,B})$.
Hence this is known as the ``earth-mover's" metric.

\begin{lemma}
\label{lemProbDist}
Let {\bf P} be a bounded region; let $W^{\psi}$ be a Wasserstein 
metric; let $\zeta$ and $\theta$ be probability distributions that are zero 
outside {\bf P}. Let $p =$diameter({\bf P}). Let
$m= \int_{{\bf x} \in \bf P} \max(0,\zeta(x)-\theta(x)) \: dx$
Then $W^{\psi}(\zeta,\theta) \leq \psi^{-1}(m \cdot \psi(p))$.
\end{lemma}

{\bf Informal proof:} The amount of ``dirt'' that has to be moved in turning
$\zeta$ into $\theta$ is \\
$ \int_{{\bf x} \in \bf P} \max(0,\zeta({\bf x})-\theta({\bf x})| \: d{\bf x}$.
The distance that any piece of dirt can be moved
is at most $p$.  So for any $\gamma$ that turns $\theta$ into $\phi$,
$I^{\psi}(\gamma,{\bf P}) \leq
m \cdot \psi(p)$.
Then $W^{\psi}(\psi,\theta) \leq 
\psi^{-1}(I^{\psi}(\gamma,{\bf P}) = 
\psi^{-1}(m \cdot \psi(p))$.

\begin{lemma}
\label{lemTopoDHFinerWasserstein2}
Let ${\bf P,Q}$ be regions. Let $p$=diameter({\bf P}), 
$h=H({\bf P,Q})$, and $a=V({\bf P,Q})$.
Assume that $a < v({\bf P})/2$ and that $h < p/2$. 
Let $\psi$ be a Mulholland function. Then
$W^{\psi}({\bf P,Q}) \leq \psi^{-1}(4a\psi(p)/v({\bf P}))$.
\end{lemma}

{\bf Proof:}
Let $\zeta = U_{\bf P}$ and $\theta = U_{\bf Q}$.
Let ${\bf R} = {\bf P} \cup {\bf Q}$; thus $\zeta$ and $\theta$ are
zero outside $\bf R$. 

Note that $v({\bf P})+a \geq v({\bf Q}) \geq v({\bf P})-a  \geq v({\bf P})/2$
\\
so $|1/(v({\bf P}) -1/v({\bf Q}))| = 
|v({\bf Q})-v({\bf P})|/(v({\bf P})v({\bf Q}))
\leq a/2v^{2}({\bf P})$.

\[
\int_{{\bf x} \in \bf R} \max(\zeta({\bf x})-\theta({\bf x}),0) \: dx =
\int_{{\bf x} \in {\bf P} \cap {\bf Q}} 
\max(\zeta({\bf x})-\theta({\bf x}),0) \: dx +
\int_{{\bf x} \in {\bf S}({\bf P},{\bf Q})} 
\max(\zeta({\bf x})-\theta({\bf x}),0) \: dx \]

But in the first integral in the sum, 
the volume of the region of integration is
at most $v({\bf P})$ and the integrand is at most $|1/v({\bf Q})-1/v({\bf P})|$
so the value of the integral is at most
$2a/v({\bf P})$.

In the second integral, the volume of integration is ${\bf S}({\bf P,Q})$
and the integrand is at most $1/\min(v({\bf P}),v({\bf Q}))$ so
value of the integral is at most $2a/v({\bf P})$.

Thus 
\[ 
\int_{{\bf x} \in \bf R} \max(\zeta({\bf x})-\theta({\bf x}),0) \: dx 
\leq 4a/v({\bf P}) \]

Using lemma~\ref{lemProbDist} it follows that 
$W^{\psi}({\bf P,Q}) \leq \psi^{-1}(4a\psi(p)/v({\bf P}))$.

\begin{theorem}
\label{thmTopoDHFinerWasserstein}
For any Mulholland function $\psi$, the topology generated
by Wasserstein distance $\Topo_{W^{\psi}}$
is coarser over $\Reg$ than the topology generated by
the dual-Hausdorff distance $\Topo_{H^{d}}$ 
\end{theorem}

{\bf Proof:}
Choose region {\bf P} and $\epsilon > 0$. Let $p$=diameter({\bf P}).
Let $b = \psi(\epsilon) v({\bf P})/4\psi(p)$. 
Using theorem~\ref{thmTopoDHFinerVolume}, choose $\delta_{1}$ such that,
such that, for all regions $\bf Q$, if $H^{d} < \delta_{1}$ then 
$V({\bf P,Q}) < b$. Let $\delta = \min(\delta_{1},p/2)$. Then
by lemma~\ref{lemTopoDHFinerWasserstein2} it follows
that $W^{\psi}({\bf P,Q}) < \epsilon$.

\begin{corollary}
\label{corWassersteinMorphing}
For any Mulholland function $\psi$, the Wasserstein distance $W^{\psi}$
supports continuous morphing over $\Reg$.
\end{corollary}

{\bf Proof:} Immediate from theorems~\ref{thmTopoDHFinerWasserstein} and
\ref{thmHausdorffContMorph}.

\begin{theorem}
\label{thmWassersteinSep}
For any Mulholland function $\psi$, the Wasserstein distance $W^{\psi}$
satisfies the region separation condition over $\Reg$.
\end{theorem}

{\bf Proof:}

{\bf Part 1:} Let ${\bf P, Z}$ be regions such that $d({\bf P,Z}) > 0$.
Let $c=d({\bf P,Z})/2$.
Let ${\bf Q}=\mbox{dilate}({\bf P},c)$. Let {\bf Y} be any 
superset of {\bf Z}. The part of {\bf Y} that is more than $c$ from
$\bf P$ includes at least $\bf Z$; the part {\bf Y} that is less than $c$
from $\bf P$ is a subset of $\bf Q$.
So the fraction of {\bf Y} that is more than 
$c$ from {\bf P} is at least $v({\bf Z})/(v({\bf Z})+v({\bf Q})$. 
So, for any uniform function $\gamma$ from {\bf P} to {\bf Y},
$I^{\psi}({\bf P},\gamma) \geq 
(v({\bf Z})/(v({\bf Z})+v({\bf Q})) \cdot \psi(c)$, so there is a positive
lower bound on $W^{\alpha}({\bf P,Y})$. 

The proof of Part 2 is analogous.

\section{The topology of the space of bounded convex open regions}
\label{secConvex}

We show that there is a unique well-behaved topology over the space of 
convex regions. Since all of the metric topologies we consider are well-behaved
over that space, it follows that they all generate the same topology.

Shephard and Webster (1995) demonstrated that the Hausdorff metric and
the symmetric-difference metric generate identical topologies over the
space of convex regions; that two further metrics, which they named the
``difference body metric'' and the ``homogeneous symmetric difference''
likewise generate the same topology. The latter two results are subsumed
in theorem~\ref{thmConvex} below, though we do not prove that here.
Groemer (2000) gives strong bounds between the relative size of
the Hausdorff distance and the symmetric-difference distance between two convex
regions.
 
\begin{lemma}
\label{lemIncreasing}
Let $\bf A$ be a bounded
open, convex region in $\Euc^{n}$. Let ${\bf p} \in \bf A$, and
let ${\bf q} \in @{\bf A}$. 
For $t \geq 0$, let ${\bf w}(t) = {\bf q}+t({\bf q}-{\bf p})$.
Then, for 
$t \geq 0$,
the function $f(t) = d({\bf w}(t), @{\bf A})$ is an increasing function of 
$t$.
\end{lemma}

\begin{figure}
\begin{center}
\includegraphics[width=4in]{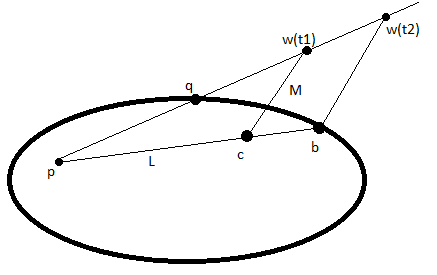}
\end{center}
\caption{Proof of lemma~\ref{lemIncreasing}}
\label{figIncreasing}
\end{figure}

{\bf Proof:} (Figure~\ref{figIncreasing}).
Let $0 < t_{1} < t_{2}$.
Let $\bf b$ be the point on $@{\bf A}$ closest to ${\bf w}(t_{2})$.  
Let $L$ be the line from {\bf p} to $\bf b$.
Since $\bf A$ is convex, the portion of $L$ between
$\bf b$ and ${\bf p}$ is entirely in $\bf A$. Let $M$ be the line through 
${\bf w}(t_{1})$ parallel to the line ${\bf bw}(t_{2})$ and let $\bf c$
be the intersection of $L$ and $M$.  Then the triangle 
$\bigtriangleup {\bf q}, {\bf w}(t_{1}),{\bf c}$ is similar to the triangle
$\bigtriangleup {\bf q}, {\bf w}(t_{2}),{\bf b}$ and lies inside it. Hence
\[ f(t_{1}) = d({\bf w}(t_{1}),{\bf A}) \leq d({\bf w}(t_{1}),{\bf c}) < 
d({\bf w}(t_{2}),{\bf b}) = f(t_{2}) \]

\begin{lemma}
\label{lemBall}
Let $\bf P$ and $\bf Q$ be bounded, convex, open sets, and let {\bf o}
be a point in {\bf P}.
Let $h=H({\bf P,Q})$ and $r=\mbox{radius}({\bf P,o})$. If $h < r$
then ${\bf B}({\bf o},r-h) \subset {\bf Q}$.
\end{lemma}

{\bf Proof:} For convenience, take $\vec{0}={\bf o}$. 
Let $\vec{x}$ be a point in ${\bf B}(\vec{0},r) \setminus {\bf Q}$
(If there is no such point, the conclusion is trivial.)
Then there is a hyperplane $\bf X$ through $\vec{x}$ such that ${\bf Q}$ lies
on one side of $\bf X$. 
Let ${\bf C}$ be the intersection of {\bf X} with 
${\bf B}(\vec{0},r)$.  ($\bf C$ is an $n-1$-dimensional solid circular
disk). Let 
$\vec{c}$ be the center of ${\bf C}$; thus $\vec{c}$ is the closest point
to $\vec{0}$ on $\bf C$, so $|\vec{c}| \leq |\vec{x}|$. 

{\bf Q} must lie in the side of {\bf X} that contains $\vec{0}$; if it
lies on the far side of {\bf X}, then its distance from the point in
${\bf B}(\vec{0},r)$ opposite $\vec{c}$ would be greater than $r$,
which is impossible.

\begin{figure}
\begin{center}
\includegraphics[width=4in]{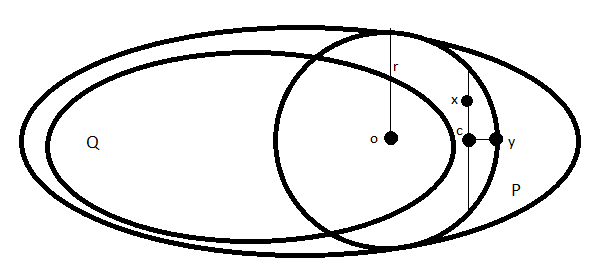}
\end{center}
\caption{Lemma~\ref{lemBall}}
\label{figLemBall}
\end{figure}

Let
$\vec{y} = r \cdot \vec{c}/|\vec{c}|$. Then $\vec{c}$ is the closest point
on ${\bf B}(\vec{0},r)$ to $\vec{y}$. In particular $d(\vec{y},\vec{c}) \leq 
d(\vec{y},{\bf Q}) \leq h$. But $d(\vec{y},\vec{c}) = r-|\vec{c}| \geq 
r-|\vec{x}|$ so $|\vec{x}| \geq r-h$, so $\vec{x} \not \in {\bf B}(\vec{0},r-h)$.
(Figure~\ref{figLemBall})
$\QED$

\begin{figure}
\begin{center}
\includegraphics[width=4in]{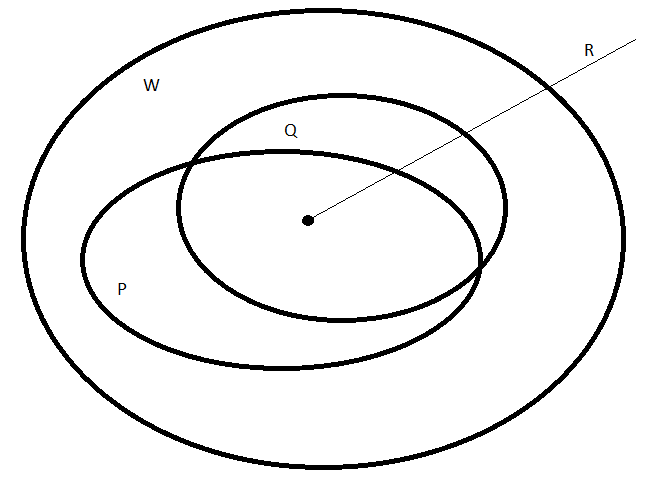}
\end{center}
\caption{The standard morphing}
\label{figStandardMorphing}
\end{figure}

\begin{definition}
\label{defStandardMorphing}
Let ${\bf P,Q,W}$ be open convex bounded regions
such that ${\bf P} \cap {\bf Q} \neq \emptyset$, 
$\bar{\bf P} \subset {\bf W}$ and 
$\bar{\bf Q} \subset {\bf W}$. That is, {\bf P} and {\bf Q} overlap, and
{\bf W} contains them both, with some separation between ${\bf P}
\cup {\bf Q}$ and the outside of {\bf W} (figure~\ref{figStandardMorphing}).

Let {\bf o} be a point in ${\bf P} \cap {\bf Q}$. 

For convenience, let $\vec{0}=\bf o$ and 
$\vec{x} = {\bf x}-{\bf o}$.
For any unit vector 
$\hat{v}$, let ${\bf R}(\hat{v})$ be the ray 
$\{ t\hat{v} \: | \: t \in (0, \infty) \}$. 
Let $\vec{p}(\hat{v})$, $\vec{q}(\hat{v})$,  $\vec{w}(\hat{v})$
be the intersections 
of ${\bf R}(\hat{v})$ with $@{\bf P}$, $@\bf Q$, and $@{\bf W}$ respectively.
Since {\bf P}, {\bf Q} and {\bf W} are convex, it is immediate that 
$\vec{p}(\hat{v})$ and $\vec{q}(\hat{v})$  and $\vec{w}(\hat{x})$ 
are uniquely defined (in
any direction $\hat{v}$ there is
only one such intersection for each) and are
continuous functions of $\hat{v}$.

The {\em standard morphing of
{\bf P} into {\bf Q} within {\bf W} centered at {\bf o}, denoted 
$\Gamma_{\bf P,Q,W,o} : [0,1] \times \Euc^{n} \mapsto \Euc^{n}$} is
defined as the following function:
\begin{quote}
For all $t \in [0,1]$, $\Gamma_{\bf P,Q,W,o}(t,\vec{0}) =\vec{0}$.

For $\vec{x} \neq \vec{0}$, let $\hat{x} = \vec{x}/|\vec{x}|$. 
To simplify the expression, fix a direction of $\hat{x}$, and let
$x = |\vec{x}|$.
$p=|\vec{p}(\hat{x})|$, $q=|\vec{q}(\hat{x})|$, and
$w=|\vec{w}(\hat{x})|$. 
Then, for any $\vec{x}$ in the ray ${\bf R}(\hat{x})$,
\begin{itemize}
\item If $x \leq p$, then
$\Gamma(t,\vec{x}) = ((1-t)x + t(xq/p)) \cdot \hat{x}$.

\item If $p < x < w$, then
$\Gamma(t,\vec{x}) = ((1-t)x+t(q+(w-q)(x-p)/(w-p))) \cdot \hat{x}$.

\item If $w \leq x$, then $\Gamma(t,\vec{x}) = \vec{x}$.
\end{itemize}
\end{quote}
\end{definition}

Thus, each ray ${\bf R}(\hat{x})$ is divided into three parts: the part inside
{\bf P}, the part between part {\bf P} and {\bf W}, and the part outside 
{\bf W}. $\Gamma$ is a transformation, piecewise
bilinear in both $t$ and $x$, which
transforms the first part into the part of the ray inside $\bf Q$, the
second part into the part of the ray between $\bf Q$ and $\bf W$, and 
is the identity outside $\bf W$.

\begin{lemma}
\label{lemExistsContMorph}
Let $\bf P, Q, W, o$ be as in definition~\ref{defStandardMorphing}.
Let $h = H({\bf P,Q})$, $r$ = radius({\bf P,o}), and
$a=\mbox{diameter}({\bf P})$.
If $r > h$, then the standard morphing $\Gamma_{\bf P,Q,W,o}$ has
the following properties:
\begin{itemize}
\item[a.] $\Gamma$ is a continuous morphing.
\item[b.] for all ${\bf x} \in \Euc^{n}$, $\Gamma(0,{\bf x})={\bf x}$.
\item[c.] for all $t \in [0,1]$ and ${\bf x} \not \in {\bf W}$,
$\Gamma(t,{\bf x})={\bf x}$.
\item[d.] $\Gamma(1,{\bf P})={\bf Q}$;
\item[e.] for all $t \in [0,1]$, 
$H(\Gamma(t,{\bf P}),{\bf P}) \leq h$; and
\item[f.] for all ${\bf x} \in \Euc^{n}$ and $t \in [0,1]$,
$d(\Gamma(t,{\bf x}),{\bf x}) \leq d({\bf x,o}) \cdot ah/(r-h)$.
\end{itemize}
\end{lemma}

{\bf Proof:}  

Properties (a), (b), and (c) are immediate by construction.

Let $\vec{0}=\bf o$;
$\vec{x} = {\bf x}-{\bf o}$ and define $\hat{x}$, 
$\vec{p}(\hat{x})$, and 
$\vec{q}(\hat{x})$ as in definition~\ref{defStandardMorphing}.

For (d):
for any point $\vec{p}(\hat{v}) \in @{\bf P}$, 
$\Gamma(0,\vec{p}(\hat{v})) = \vec{p}(\hat{v})$ and
$\Gamma(1,\vec{p}(\hat{v})) = \vec{q}(\hat{v})$.
Since {\bf P} and {\bf Q} are convex, it follows that 
$\Gamma(1,@{\bf P})=@{\bf Q}$ and therefore
$\Gamma(1,{\bf P})={\bf Q}$.

Condition (e) of the lemma asserts that,
for all $t$, $H({\bf P},\Gamma(t,{\bf P})) \leq
H({\bf P},{\bf Q})$; that is
for all $\vec{x} \in \Gamma(t,{\bf P})$, $d(\vec{x},{\bf P}) \leq h$
and
for all $\vec{x} \in {\bf P}$, $d(\vec{x},\Gamma(t,{\bf P})) \leq h$

To prove this,
let $\vec{x}$  be a point in $\Gamma(t,{\bf P})$, 
and let $\hat{x}=\vec{x}/|\vec{x}|$.
Then the points $\vec{0}, \vec{x}, \vec{p}(\hat{x})$, 
and $\vec{q}(\hat{x})$ are
collinear. 
If $|\vec{x}| < |\vec{p}(\hat{x})|$ then $\vec{x} \in {\bf P}$, so
$d(\vec{x},{\bf P})=0$. 
If $|\vec{x}| \geq |\vec{p}(\hat{x})|$ 
then $|\vec{q}(\hat{x})| > |\vec{p}(\hat{x})|$ and
$\vec{x}$ is on the line between $\vec{p}(\hat{x})$ and $\vec{q}(\hat{x})$ 
so, by lemma~\ref{lemIncreasing}, 
$d(\vec{x},{\bf P}) \leq d(\vec{q}(\hat{x}),{\bf P}) \leq H({\bf Q},{\bf P})$.

Now let $\vec{x}$  be a point in $\bf P$, and let $\hat{x}=\vec{x}/|\vec{x}|$.
If $\vec{x} \in \Gamma(t,{\bf P})$ then $d(\vec{x},f(t,{\bf P}))=0$.
If $\vec{x} \not \in \Gamma(t,{\bf P})$ then $\vec{x}$ must be on the line through
$\vec{q}(\hat{x})$ and $\vec{p}(\hat{x})$ with
$|\vec{q}(\hat{x})| < |\vec{x}| < \vec{p}(\hat{x})$. By lemma~\ref{lemIncreasing}
$d(\vec{x},{\bf Q}) \leq d(\vec{p}(\hat{x}),{\bf Q}) \leq H({\bf Q},{\bf P})$.

Condition (f) of the lemma asserts that 
for all $\vec{x} \in \Euc^{n}$ and $t \in [0,1]$,
$d(\Gamma(t,\vec{x}),\vec{x}) \leq ph/(r-h)$. 
By construction,
the point on the ray $\{ t\hat{x} | t > 0 \}$ that is moved furthest
is $\vec{p}(\hat{x})$, so it suffices to prove the inequality for that
point.

Since
\[ \Gamma(t, \vec{x}) = \vec{x} \cdot (1 + 
\frac{t \cdot (|\vec{q}(\hat{x}|-|\vec{p}(\hat{x}|}{|\vec{p}(\hat{x})|}) \]
we have 
\[ d(\Gamma(t,\vec{x}),\vec{x}) = |\vec{x}| \cdot  t
\frac{\mbox{abs}(|\vec{q}(\hat{x}|-|\vec{p}(\hat{x}|)}{|\vec{p}(\hat{x})|} \]
Our goal,  then, is to bound the above fraction as a function of $r$ and $h$.
For convenience since $\hat{x}$ will be fixed, we will drop the argument and
just write $\vec{p}$ and $\vec{q}$.

\begin{figure}
\begin{center}
\includegraphics[width=4in]{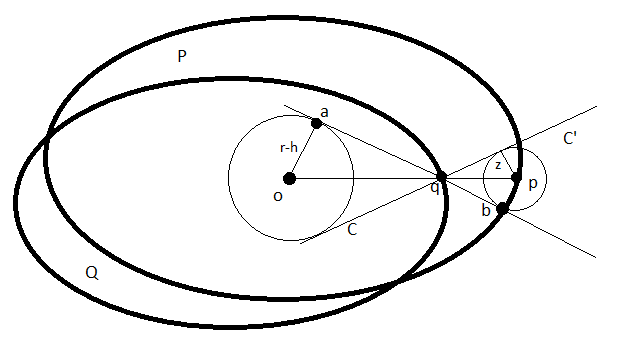}
\caption{Proof of lemma~\ref{lemExistsContMorph}}
\end{center}
\label{figCone}
\end{figure}

Consider first the case where $|\vec{q}| < |\vec{p}|$. 
The ray ${\bf R} = \{ t \hat{x} | t \in (0, \infty) \}$ is thus 
divided into three parts: the segment  from $\vec{0}$ to
$\vec{q}$ is in both {\bf Q} and {\bf P}; the
segment from $\vec{q}$ to $\vec{p}$ 
is in ${\bf P}$ but not $\bf Q$; and the segment
past $\vec{p}$ is in neither. By 
lemma~\ref{lemBall}, the ball ${\bf B}(\vec{0},r-h) \subset {\bf Q}$.
Construct the cone $\bf C$ with apex $\vec{q}$ that is tangent to
${\bf B}(\vec{0},r-h)$ (figure~\ref{figCone}). 
Since $\bf Q$ is convex, ${\bf C} \subset {\bf Q}$.
Let $\bf C'$ be the reflection of {\bf C} through $\vec{q}$.
Then $\bf C'$ must be disjoint from $\bf Q$. 
(For any point $\vec{w} \in {\bf C'}$
there are points $\vec{v}$ on the part of the ray {\bf R} past $\vec{q}$ 
and $\vec{u} \in {\bf C}$ such that $\vec{u},\vec{v}, \vec{w}$ are
collinear in that order; since $\bf Q$ is convex, $\vec{u} \in {\bf Q}$ and
$\vec{v} \not \in \bf Q$, it follows that $\vec{w} \not \in {\bf Q}$.)

Construct the sphere centered at $\vec{p}$ tangent to ${\bf C}'$. Let 
$z$ be the radius of the sphere. Since ${\bf p} \in @{\bf P}$ and
the sphere is disjoint from $\bf Q$, we have $h \geq z$. 

Now let 
$\vec{a}$ be a point in $\bar{\bf B}(\vec{0},(r-h)) \cap {\bf C}$
and let
$\vec{b}$ be a point in $\bar{\bf B}(\vec{p},z) \cap {\bf C}'$ such that
$\vec{a}, \vec{q}, \vec{b}$ are collinear. Then the triangles
$\bigtriangleup \vec{0},\vec{a},\vec{q}$ and
$\bigtriangleup \vec{p},\vec{b},\vec{q}$ are similar right triangles.
So $d(\vec{0},\vec{a})/d(\vec{0},\vec{q}) = (r-h)/|\vec{p}| =
d(\vec{p},\vec{b})/d(\vec{p},\vec{q}) = z/(|\vec{p}| - |\vec{q}|)$.

Combining these and rearranging we get 
$(|\vec{p}|-|\vec{q}|)/|\vec{p}| \leq h/(r-h)$.

In the case where $|\vec{p}| < |\vec{q}|$, the analysis is exactly
analogous, except that in that case you get the tighter bound
$(|\vec{p}|-|\vec{q}|)/|\vec{p}| \leq h/r$.

$\QED$

\begin{corollary}
\label{corMapping}
Let ${\bf P,Q}$ be convex regions such that ${\bf P} \cap {\bf Q} \neq
\emptyset$. 
Let $h = H({\bf P,Q})$, $r$ = radius({\bf P,o}), and
$a=\mbox{diameter}({\bf P})$.
Then there is a homeomorphism $g$ of $\Euc^{n}$ to itself such 
that $g({\bf P}) = {\bf Q}$ and, 
for all ${\bf x} \in {\bf P}$, $d({\bf x},g({\bf x})) \leq ah/(r-h)$.
\end{corollary}

{\bf Proof:} Find a convex region ${\bf w} \supset {\bf P} \cup {\bf Q}$
and choose a point ${\bf o} \in {\bf P} \cap {\bf Q}$. Then by
lemma~\ref{lemExistsContMorph} the function
$\Gamma_{\bf P,Q,W,o}(1, \cdot)$ satisfies the condition of the corollary.

It seems likely that this bound can be substantially tightened using
a different morphing and in 
particular that the dependence on diameter({\bf P}) can be eliminated.
But for the purposes of our analysis, this will suffice.

\begin{lemma}
\label{lemContMorph}
Let $\Topo$ be a topology over $\Reg$ that supports
continuous morphing. Then, restricted to $\Convex$, $\Topo_{H}$, 
the topology induced by the Hausdorff metric, is at least as fine as $\Topo$.
\end{lemma}

{\bf Proof} of the contrapositive:  Suppose that $\Topo_{H}$ is not
a refinement of $\Topo$. Then there exists a region ${\bf P} \in \Convex$
and a sequence of regions ${\bf Q}_{1}, {\bf Q}_{2} \ldots \in \Convex$
that converges to ${\bf P}$ in $\Topo_{H}$ but not in $\Topo$. 
Let $r = \mbox{radius}({\bf P}) > 0$.
Let $\epsilon_{i} = H({\bf Q}_{i},{\bf P})$; thus $\lim_{i \goesto \infty}
\epsilon_{i} = 0.$ By renumbering we can assume that 
$\epsilon_{i} < r/2$ for all $i$.

We are going to use lemma~\ref{lemExistsContMorph} to interpolate a continuous
morphing $\phi$ that passes through the 
regions ${\bf Q}_{1}, {\bf Q}_{2}, {\bf Q}_{3} \ldots {\bf P}$ 
at times 1, 1/2, 1/3 \ldots 0. 

Fix a center point ${\bf o} \in {\bf P}$. 
By lemma~\ref{lemBall}, ${\bf B}({\bf o},r/2) \subset
{\bf B}({\bf o},r-H({\bf Q}_{i},{\bf P})) \subset {\bf Q}_{i}$.

Let 
$q = 1 + \mbox{diameter}({\bf P})+\max_{i}H({\bf Q}_{i},{\bf P})$;
then it is easily shown that the sphere ${\bf R}={\bf B}({\bf o},q)$
contains $\bar{\bf P}$ and $\bar{\bf Q}_{i}$ for all $i$.

Define the function $f_{k} = \Gamma_{{\bf Q}_{k},{\bf Q}_{k+1}, {\bf R},
{\bf o}}$ as in definition~\ref{defStandardMorphing}.
By lemma~\ref{lemExistsContMorph}, 
$f_{k}(t,{\bf x})$ is a continuous morphing,
$f_{k}(0,\cdot)$ is the identity, and 
$f_{k}(1,{\bf Q}_{i}) ={\bf Q}_{i+1}$. \\
Define the function $g_{k}(t,{\bf x}) = f_{k}(k+1-k(k+1)t, {\bf x})$; thus
$g_{k}(1/k,{\bf x}) = f_{k}(0,{\bf x})$ and
$g_{k}(1/(k+1),{\bf x}) = f_{k}(1,{\bf x})$.

Now define the function $\phi: \Real \times \Euc^{n} \mapsto \Euc^{n}$ as
follows:
\begin{itemize}
\item Construct $f_{0}$ to satisfy lemma~\ref{lemExistsContMorph} for 
{\bf P} and {\bf Q}. For $t \geq 1$, define $\phi(t,{\bf x)} = 
f_{0}(1,{\bf x})$.
\item For $k=1,2,3 \ldots$, for $t \in [1/(k+1),1/k)$ define
$\phi(t,{\bf x}) = g_{k}(t,(\phi(1/k,{\bf x}))$
\item for $t \leq 0$, $\phi(t,\cdot)$ is the identity function on
$\Euc^{n}$
\end{itemize}

Note that 
$\phi(1,{\bf P}) = f_{0}(1,{\bf P}) = {\bf Q}_{0}$. \\
$\phi(1/2,{\bf P}) = g_{1}(1/2,\phi(1,{\bf P})) = 
f_{1}(1,Q_{0}) = {\bf Q}_{1}$. \\
$\phi(1/3,{\bf P}) = g_{2}(1/3,\phi(1/2,{\bf P})) = 
f_{2}(1,Q_{1}) = {\bf Q}_{2}$. \\
and in general $\phi(1/k,{\bf P}) = {\bf Q}_{k}$.

To show that $\phi$ is continuous: Spatial continuity is immediate by 
construction.
Temporal continuity between times of the form
$1/k$ is guaranteed by the continuity of $f_{k}$. Continuity at times
of the form $1/k$ follows from the fact that $\phi(t,\cdot)$ consists
in expansion along rays emanating from a fixed center point $\vec{0}$
and that the limit at time $t=1/k$, both from above and below, of the
amount of expansion at point $\vec{x}$ is 
$|\vec{q}_{k}(\hat{x})| / |\vec{p}(\hat{x})|$, in the notation of 
lemma~\ref{lemExistsContMorph}, where $\vec{q}_{k}(\hat{x})$ is the intersection
of ${\bf Q}_{k}$ with the ray $\{ t \cdot \hat{x} \: | \: t > 0 \}$.

The continuity of $\phi$ at time $t=0$, which is, of course, the critical
point,
is guaranteed by the facts that, by lemma~\ref{lemExistsContMorph},
for all $t \in [1/(k+1),1/k]$, 
$d(\phi(t,\vec{x}),\phi(1/(k+1),\vec{x}) \leq
2H({\bf Q}_{k},{\bf Q}_{k+1})/r$, and that 
$d(\phi(1/(k+1),\vec{x}),\phi(0,\vec{x}) \leq
2H({\bf Q}_{k},{\bf P})/r$, and by assumption, both of these Hausdorff distances
go to zero as $k \goesto \infty$.

$\QED$

\begin{lemma}
\label{lemSequence}
Let $\bf P$ be a bounded open region and let ${\bf Q}_{1}, 
{\bf Q}_{2} \ldots $ be an 
infinite sequence of convex, open regions. Then one of three things is true.
\begin{itemize}
\item[1.] $\lim_{i \goesto \infty} H({\bf P,Q}_{i}) = 0$.
\item[2.] There is a region $\bf Z$ such that 
${\bf Z} \subset {\bf P}$ and,
for infinitely many ${\bf Q}_{i}$, ${\bf Z }\cap {\bf Q}_{i} = \emptyset$.
\item[3.] There is a region $\bf Z$ such that 
${\bf Z} \cap {\bf P} = \emptyset$ and,
for infinitely many ${\bf Q}_{i}$, ${\bf Z} \subset {\bf Q}_{i}$.
\end{itemize}
\end{lemma}

{\bf Proof:} If condition 1 does not hold, then there exists $c >0$ such
that either
(a) $H^{1}({\bf P,Q}_{i}) > c$ for infinitely many $i$, or
(b) $H^{1}({\bf Q}_{i},{\bf P}) > c$ for infinitely many $i$.

Suppose that (a) holds. For each such ${\bf Q}_{i}$, there is a point
${\bf p}_{i} \in {\bf P}$ such that $d({\bf p}_{i},{\bf Q}_{i}) > c$. 
These ${\bf p}_{i}$ must have a cluster point ${\bf p}$ in the 
closure of {\bf P}. Choose $\epsilon$ so that $0 < \epsilon < c$,
and let the infinite set of indices $I =\{ i \: | \:  d({\bf p}_{i},{\bf p}) <
\epsilon \}$. Then for $i \in I$, $d({\bf p},{\bf Q}_{i}) > c-\epsilon$.
Therefore condition 2 of the lemma is satisfied for 
${\bf Z} = {\bf P} \cap {\bf B}({\bf p},c-\epsilon)$.

Suppose that conditions 1 and 2 and (a) do not hold but (b) holds. 
Since ${\bf P}$
is open,
there exists an open ball ${\bf B}({\bf o},r) \subset {\bf P}$.
Let $0 < \epsilon  < r$.
Since (a) does not hold, $H({\bf P},{\bf Q}_{i}) < \epsilon$
for all but finitely many $i$. Ignore the $i$ where it does not happen.
By lemma~\ref{lemBall}, $B({\bf p},r-\epsilon) \subset {\bf Q}_{i}$. 
Let $r'=\min(c,r-\epsilon)$. 

Since ${\bf P}$ is bounded, let $s$ be such that 
${\bf P} \subset {\bf B}({\bf o},s)$. 

Since case (b) holds, for each ${\bf Q}_{i}$ there is a point
${\bf q}_{i} \in {\bf Q}_{i}$ such that $d({\bf q}_{i},{\bf P})  > c$.

Let ${\bf H}_{i}$ be the convex hull of 
${\bf B}({\bf o},r') \cup {\bf B}({\bf q}_{i},r')$. Thus ${\bf H}_{i}$ is
a right spherical cylinder with spherical caps whose axis is the line
from {\bf o} to ${\bf q}_{i}$. Since 
${\bf B}({\bf o},r') \subset {\bf Q}_{i}$,
${\bf B}({\bf q}_{i},r') \subset {\bf Q}_{i}$, and ${\bf Q}_{i}$ is convex,
${\bf H}_{i} \subset {\bf Q}_{i}$. 

\begin{figure}
\begin{center}
\includegraphics[width=4in]{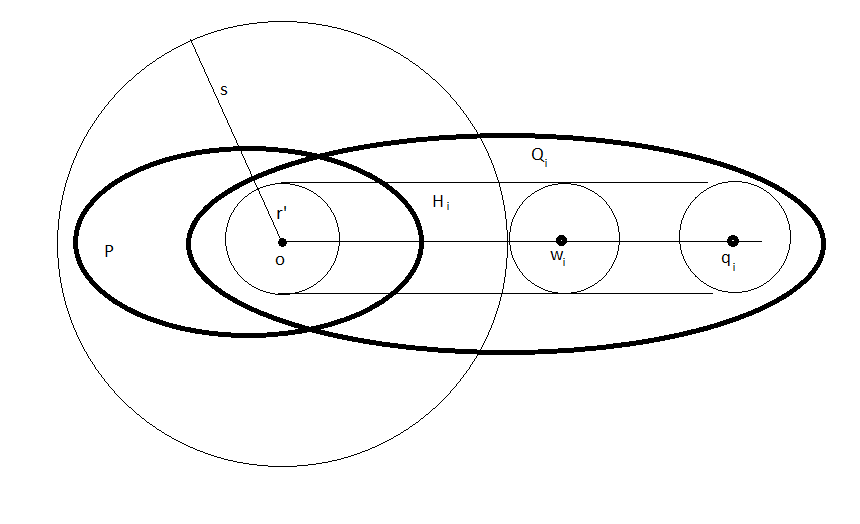}
\end{center}
\caption{Lemma~\ref{lemSequence}: Condition 3}
\label{figLemSequence}
\end{figure}

Let ${\bf w}_{i} = {\bf o}+\min(1,(s+c)/d({\bf q}_{i},{\bf o})) \cdot 
({\bf q}_{i}-{\bf o})$; that is ${\bf w}_{i}$ is either ${\bf q}_{i}$,
if ${\bf q}_{i}$ is less than distance $s+r'$ from  {\bf o} or is the
point on the line from ${\bf o}$ to ${\bf q}_{i}$ at distance $s+c$
from ${\bf o}$. In either case, 
${\bf Z} = {\bf B}({\bf w}_{i},r')$ is disjoint from
$\bf P$ and is a subset of ${\bf H}_{i}$ and therefore of ${\bf Q}_{i}$ 
(figure~\ref{figLemSequence}).

Since all the ${\bf w}_{i}$ lie in the bounded region 
$\bar{\bf B}({\bf o},s+r')$,
they have a cluster point $\bf w$. Thus, for any $t < r'$, 
${\bf B}({\bf w},t)$
is a subset of infinitely many ${\bf Q}_{i}$ and is disjoint from {\bf P}.

$\QED$

\begin{lemma}
\label{lemConvexAtLeastAsFine}
Let  $\mu$ be a metric on $\Reg$ such that the topology $\Topo_{\mu}$ 
satisfies the region separation condition. 
Then over the space of convex open regions,
$\Topo_{\mu}$ is at least as fine as 
$\Topo_{H}$, the
topology of the Hausdorff metric.  
\end{lemma}

{\bf Proof by contradiction:} Suppose that $\Topo_{\mu}$ is not at least
as fine as $\Topo_{H}$. Then there exists $\epsilon > 0$ and a region $\bf P$
such that the ball in the Hausdorff-metric topology 
${\mathcal B}_{H}({\bf P},\epsilon)$
is not contained in any ball in the $\mu$ topology. Thus, there is a sequence
of regions ${\bf Q}_{1}, {\bf Q}_{2} \ldots$ such that 
$\mu({\bf Q}_{i},{\bf P}) < 1/i$ but
$H({\bf Q}_{i},{\bf P}) \geq \epsilon$ for all $i$. By lemma~\ref{lemSequence}
either (a)
there exists a region ${\bf Z} \subset \bf P$ such that 
$\bf Z$ is disjoint from ${\bf Q}_{i}$ for infinitely many ${\bf Q}_{i}$; or
(b) there exists a region ${\bf Z}$ disjoint from $\bf P$ such that 
${\bf Z} \subset {\bf Q}_{i}$ for infinitely many ${\bf Q}_{i}$. 

Let ${\mathcal U} \in \Topo_{\mu}$ satisfy the conditions of 
definition~\ref{defSeparates}. Then by that definition, infinitely many 
${\bf Q}_{i}$ are not in $\mathcal U$; but that contradicts their
construction above.

\begin{theorem}
\label{thmConvex}
Let $\Topo_{\mu}$ be a well-behaved metric topology.
Over the space $\mathcal C$ of convex open regions, $\Topo_{\mu}$ is
equal to $\Topo_{H}$, the topology of the Hausdorff metric.
\end{theorem}

{\bf Proof:} This is just the combinations of lemmas~\ref{lemContMorph} 
and \ref{lemConvexAtLeastAsFine}.

\begin{corollary}
\label{corIdentical}
Over the space $\mathcal C$ of convex open regions, the metrics 
$M,H,H^{d},V$ and $W^{\psi}$ all generate the identical topology.
\end{corollary}

{\bf Proof:} Immediate from theorem~\ref{thmConvex} together with 
theorems~\ref{thmMContinuous}, \ref{thmHausdorffContMorph}, 
\ref{thmHausdorffSeparation}
\ref{thmVolumeSupportsMorphing}
\ref{thmVolumeSeparation}
\ref{thmWassersteinSep} and corollary\ref{corWassersteinMorphing}.

\section{The space of two separated convex regions}
\label{secTwoConvex}
We now turn to, arguably, the next simplest class of regions: those that
consist either of a single convex region or are the union of two separated 
convex regions. As we shall see, our metrics generate many different topologies
for that space.

Let $\mathcal D^{2}$ be the set of all unions of two separated 
convex regions:
${\mathcal D}^{2} = \{ {\bf X} \cup {\bf Y} \: | \:
{\bf X,Y} \in {\mathcal C}, d({\bf X,Y}) > 0 \}$. \\
Let ${\mathcal D} = {\mathcal C} \cup {\mathcal D}^{2}$.

\subsection{Well-behaved topologies over $\mathcal D$}
We begin by establishing some properties of any well-behaved 
topology over $\mathcal D$.

Let {\bf A} be a region in $\mathcal D$ and let $\Topo$ be a well-behaved
topology over $\mathcal D$. Theorem~\ref{thmConvex} above showed that, informally,
speaking, if {\bf A} is convex, the convex regions close to {\bf A} in $\Topo$
are those that are close in the Hausdorff distance. We will show in 
that, if {\bf A} is ${\mathcal D}^{2}$, 
then small neighborhoods of {\bf A} contain no convex regions
(lemma~\ref{lemNoConvexClosetoNonConvex})
and that they contain
exactly the regions in ${\mathcal D}^{2}$ that are close in the Hausdorff
distance (theorem~\ref{thmNonConvexHausdorff}). 
The interesting question is, if {\bf A} is convex, what kinds of
regions in ${\mathcal D}^{2}$ lie in its neighborhoods?  As we will see,
there are many different possible answers, depending on the metric.

\begin{lemma}
\label{lemNonConvexFarFromConvex}
Let $\bf P$ be a region that is not convex. Then there 
exists  $\epsilon>0$ such
that, for every convex region $\bf Q$, radius({\bf S}({\bf P,Q})) $\geq$
$\epsilon$.
\end{lemma}

{\bf Proof:} Since $\bf P$ is not convex, let $\bf a,b,c$ be points such that 
{\bf b} lies on line {\bf ac}, ${\bf a,c} \in \bf P$ and ${\bf b} \not \in
\bar{\bf P}$. Let $\epsilon_{1}> 0$ be such that 
${\bf B}({\bf a},\epsilon_{1}) \subset {\bf P}$, 
${\bf B}({\bf c},\epsilon_{1}) \subset {\bf P}$,  and
${\bf B}({\bf b},\epsilon_{1})$ is disjoint from $\bar{\bf P}$. 
If both {\bf a} and {\bf c} are in {\bf Q}, then {\bf b} is in {\bf Q},
so $H({\bf P,Q}) \geq d({\bf b,Q}) \geq \epsilon$.
If {\bf a} is not in {\bf Q}, then, since {\bf Q} is convex,
some hemisphere of ${\bf B}({\bf a},\epsilon_{1})$ is not in {\bf Q}. 
This hemisphere contains a ball of radius $\epsilon_{1}/2$.
The same holds if {\bf c} is not in {\bf Q}. Therefore,
the conclusion is satisfied with $\epsilon = \epsilon_{1}/2$.

\begin{lemma}
\label{lemNoConvexClosetoNonConvex}
Let $\mu$ be either the Hausdorff metric, the symmetric difference metric,
or a Wasserstein metric.
Let $\bf P$ be a non-convex region. Then there
exists $\epsilon > 0$ such that ${\mathcal B}_{\mu}({\bf P},\epsilon)$ 
does not contain any convex regions.
\end{lemma}

{\bf Proof:} 
Immediate from lemma~\ref{lemNonConvexFarFromConvex}.

\begin{lemma}
\label{lemMatchPieces}
Let {\bf P}={\bf C} $\cup$ {\bf D}  and {\bf Q}={\bf E} $\cup$ {\bf F}. 
where {\bf C, D, E,} and {\bf F} are convex, 
$d({\bf C,D}) > 0$, and $d({\bf E,F}) > 0$.
Let $r_{C}$ and $r_{D}$ be the radii of ${\bf C}$ and {\bf D} respectively. 
Let $h = H({\bf P,Q})$. 
If $h < \min(r_{\bf C},r_{\bf D},d({\bf C,D})/2)$, then either
\begin{itemize}
\item[a.] 
radius(${\bf C} \cap {\bf E}$) $>$ $r_{C}- h$, 
$H({\bf C},{\bf E}) \leq h$,
${\bf C} \cap {\bf F} = \emptyset$, 
radius(${\bf D} \cap {\bf F}$) $>$ $r_{D}- h$,
$H({\bf D},{\bf F}) \leq h$, and
${\bf D} \cap {\bf E}=\emptyset$; or
\item[b.]
radius(${\bf D} \cap {\bf E}$) $>$ $r_{D}- h$, 
$H({\bf D},{\bf E}) \leq h$,
${\bf D} \cap {\bf F} = \emptyset$,
radius(${\bf C} \cap {\bf F}$) $>$ $r_{C}- h$, 
$H({\bf C},{\bf F}) \leq h$, and
${\bf C} \cap {\bf E}=\emptyset$
\end{itemize}
In case (a), we say that {\em {\bf E} corresponds to {\bf C}
and {\bf F} to {\bf D}}.
\end{lemma}

{\bf Proof:} Let {\bf c} be a point such that 
${\bf B}({\bf c},r_{\bf C}) \subset {\bf C}$. Since $H^{1}({\bf P,Q}) \leq h$,
there is a point ${\bf q} \in {\bf Q}$ such that $d({\bf c,q}) < h$, so 
${\bf q} \in {\bf C}$.  Since ${\bf Q} = {\bf E} \cup {\bf F}$, it follows that
${\bf q} \in {\bf E}$ or
${\bf q} \in {\bf F}$; let us say in {\bf E}.

I claim that $d({\bf D,E}) > h$. Proof by contradiction. 
Suppose there are points ${\bf d} \in {\bf D}$ and
${\bf e} \in {\bf E}$ such that $d({\bf d,e}) \leq h$.
Let ${\bf z}$ be the point in $\bar{C}$ closest to {\bf e}; then
$d({\bf e,C}) = d({\bf e,z})$. Also $d({\bf C,D}) \leq d({\bf z,d}) \leq
d({\bf z,e})+d({\bf e,d})$. By assumption of the lemma,
$2h < d({\bf C,D})$. Combining these we have $d({\bf e,C}) > h$. 

For any point {\bf x} let $\phi({\bf x}) = d({\bf x,C})-d({\bf x,D})$.
As you move on a straight line from {\bf q} to {\bf e}, the value of
$\phi$ changes from positive to negative. Let {\bf y} be a point where 
$\phi({\bf y})=0$ so $d({\bf y,D}) = d({\bf y,C})$. Again we have
inequality that $2h < d({\bf y,C}) + d({\bf y,D})$ so 
$d({\bf h,P}) = \min(d({\bf y,C}), d({\bf y,D})) > h$.
Since $H^{1}({\bf E,P}) \leq h$ that means that {\bf y} is not in
{\bf E}. But since {\bf E} is convex, and {\bf q} and {\bf e} are 
in {\bf E}, {\bf y} must be in {\bf E}. That completes the contradiction.

Since $H^{1}({\bf D,Q}) \leq  h$ and $d({\bf E},{\bf D}) > h$,
it must be that $H^{1}({\bf F,D}) \leq h$. It follows from lemma~\ref{lemBall}
that radius(${\bf F} \cap {\bf D}$) $\geq$ $r_{\bf D}-h$.

The same arguments show that $d({\bf E,D}) > h$ and 
that radius(${\bf E} \cap {\bf C}$) $\geq$ $r_{\bf C}-h$.

$\QED$
\begin{lemma}
\label{lemH1OfConvexHull}
Let {\bf P} be a convex region; let {\bf Q} be a region; and let {\bf R}
be the convex hull of ${\bf P} \cup {\bf Q}$. Then 
$H^{1}({\bf R,P}) = H^{1}({\bf Q,P})$ 
\end{lemma}

{\bf Proof:} Let {\bf r} be the point in $\bar{\bf R}$ that is furthest
from {\bf P}. There exists points ${\bf u,v} \in \bar{\bf P} \cup \bar{\bf Q}$
such that $\bf r$ is on the line {\bf uv}. Let {\bf w,x} be the points in
$\bar{\bf P}$ closest to {\bf u,v} respectively. Since {\bf P} is
convex, the line {\bf wx} is in {\bf P}. It is always the case that, given
two lines {\bf uv} and {\bf wx} and a point {\bf r} on {\bf uv}, 
$d({\bf r,wx}) \leq \max(d({\bf u,w}),d({\bf v,x}))$. (The distance squared
is a convex quadratic function, whose maximum over any interval is reached
at one of the extrema.) So we have
$H^{1}({\bf R,P}) = d({\bf r,P}) \leq d({\bf r,wx}) \leq 
\max(d({\bf u,w}),d({\bf v,x}))) \leq H^{1}({\bf Q,P})$.
The reverse inequality is trivial.

\begin{lemma}
\label{lemExistsContMorphConvex2}
(Analogous to lemma~\ref{lemExistsContMorph}). 
Let $\bf P, Q$ be regions in ${\mathcal D}^{2}$.
Let ${\bf C,D,E,F}$ be convex regions such that 
${\bf P}={\bf C} \cup {\bf D}$; 
${\bf Q}={\bf E} \cup {\bf F}$; 
{\bf E} corresponds to {\bf C} and {\bf F} corresponds to {\bf D}.
Let $h = H({\bf P,Q})$. 
Let $r = \min(\mbox{radius}({\bf C}),\mbox{radius}({\bf D}))$ and let 
$p = \max(\mbox{diameter}({\bf C}),\mbox{diameter}({\bf D}))$.
If $h < d({\bf C,D})/2$ then there exists a continuous 
morphing $f:[0,1] \times R^{n} \mapsto R^{n}$ such that:
\begin{itemize}
\item[a.] for all ${\bf x} \in \Euc^{n}$ $f(0,{\bf x})={\bf x}$.
\item[b.] $f(1,{\bf P})={\bf Q}$;
\item[c.] for all $t \in [0,1]$, 
$H(f(t,{\bf P}),{\bf P}) \leq h$; and
\item[d.] for all ${\bf x} \in \Euc^{n}$ and $t \in [0,1]$,
$d(f(t,{\bf x}),{\bf x}) \leq d({\bf x,o}) \cdot ph/(r-h)$.
\end{itemize}
\end{lemma}

{\bf Proof:} 
Let {\bf W} be the convex hull of ${\bf C} \cup {\bf E}$
and let {\bf X} be the convex hull of ${\bf D} \cup {\bf F}$.
By lemma~\ref{lemH1OfConvexHull} $H^{1}({\bf W,C}) \leq  h$ and
$H^{1}({\bf X,D}) \leq  h$.
Let $\epsilon = d({\bf C,D})-2h > 0$.
Let $\bf R$ and $\bf S$ be the expansions of {\bf W} and {\bf X} by
$\epsilon$; that is
${\bf R} = \{ {\bf r} \: | \: d({\bf r},{\bf W}) < \epsilon \}$ and`
${\bf S} = \{ {\bf r} \: | \: d({\bf r},{\bf X}) < \epsilon$.
It is easily shown that {\bf R} and {\bf S} are convex and disjoint.

Choose points ${\bf c} \in {\bf C}$, ${\bf d} \in {\bf D}$ such that
${\bf B}({\bf c},r) \subset {\bf C}$, 
${\bf B}({\bf d},r) \subset {\bf D}$. 
Clearly 
the maximal distance from {\bf c} to a point on @{\bf C} and
the maximal distance from {\bf d} to a point on @{\bf D} are
at most $p$.
 
We can use definition~\ref{defStandardMorphing} to construct functions
$\Gamma_{\bf C,E,R,c}$ and $\Gamma_{\bf D,F,S,d}$. Define $f(t,{\bf x})$ as
\[ f(t,{\bf x}) = \left\{ \begin{array}{ll}
\Gamma_{\bf C,E,R,c}(t,{\bf x}) & \mbox{if } {\bf x} \in {\bf R} \\
\Gamma_{\bf D,F,S,d}(t,{\bf x}) & \mbox{if } {\bf x} \in {\bf S} \\
{\bf x} & \mbox{otherwise}
\end{array} \right. \]

The stated properties then follow immediately from the properties of
$\Gamma$ in lemma~\ref{lemExistsContMorph}.

\begin{theorem}
\label{thmNonConvexHausdorff}
Let $\Topo_{\mu}$ be a well-behaved metric topology.
Then the restriction of
$\Topo_{\mu}$ to ${\mathcal D}^{2}$ is equal to $\Topo_{H}$, the
topology of the Hausdorff metric.
\end{theorem}

{\bf Proof:} Identical to the proof of theorem~\ref{thmConvex}, replacing
the use of lemma~\ref{lemExistsContMorph}
with lemma~\ref{lemExistsContMorphConvex2}.

Thus, in view of theorems~\ref{thmConvex} and \ref{thmNonConvexHausdorff} and  
lemma~\ref{lemNoConvexClosetoNonConvex}, if $\Topo_{\mu}$ is the Hausdroff,
the symmetric difference, or the Wasserstein metric
topology over $\mathcal D$, then
every neighborhood of a region in ${\mathcal D}^{2}$ is a set of regions, all
in ${\mathcal D}^{2}$ that are close in the Hausdorff metric; while the convex
regions in the neighborhood of a convex region are those that are close
in the Hausdorff distance. All that remains, therefore, is to characterize
the non-convex regions that lie in the neighborhood of a convex region.
We now explore how that works out in the various metrics we are studying.

\subsection{The homeomorphism-based topology in $\mathcal D$}

Over the space $\mathcal D$, the topology $\Topo_{M}$ is uninteresting.
The distance between a region in $\mathcal C$ and
a region in ${\mathcal D}^{2}$ is always infinite, so a basis for the
topology over $\mathcal D$ is 
(the open sets of the Hausdorff topology over $\mathcal C$) union
(the open sets of the Hausdorff topology over $\mathcal D$).
In other words, the question, ``What regions in ${\mathcal D}^{2}$ are close
to a convex region in $\mathcal C$?'' has the most boring possible answer:
None at all.

\subsection{The dual-Hausdorff metrics in ${\mathcal D}$}
The dual-Hausdorff metric topology is strictly coarser than the homeomorphism
metric topology over $\mathcal D$.
In particular, a history in which a growing, second, piece emerges from the 
surface of a convex region is continuous under $H^{d}$. Thus, histories 
1 and 2 are continuous under $H^{d}$ but not under $M$.

\begin{figure}
\begin{center}
\includegraphics[width=4in]{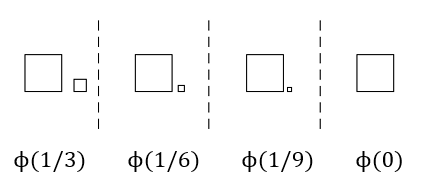}
\end{center}
\caption{History 1.0}
\label{figHistory1.0}
\end{figure}

{\bf History 1.0:} In $\Euc^{2}$
let $\phi(0) = (0,1) \times (0,1)$.  For $t > 0$, let 
$\phi(t) = (0,1) \times (0,1) \cup (1+t,1+2t) \times (0,t)$
(figure~\ref{figHistory1.0}).

\begin{figure}
\begin{center}
\includegraphics[width=4in]{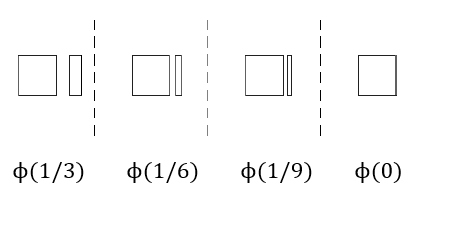}
\end{center}
\caption{History 1.1}
\label{figHistory1.1}
\end{figure}

{\bf History 1.1:} In $\Euc^{2}$
let $\phi(0) = (0,1) \times (0,1)$.  For $t > 0$, let 
$\phi(t) = (0,1) \times (0,1) \cup (1+t,1+2t) \times (0,1)$
(figure~\ref{figHistory1.1}).

It seems somewhat plausible that, for some purpose, one might consider
history 1.0 to be continuous, but not history 1.1. This can be achieved 
in $\Euc^{2}$ as follows:
Let perimeter($\bf P$) be the perimeter of region $\bf P$ (i.e. the arc length
of @{\bf P}). Define a metric $\mu$ as follows:
$\mu({\bf P,Q}) = H^{d}({\bf P,Q}) + 
\mbox{abs}(\mbox{perimeter}({\bf P})-\mbox{perimeter}({\bf Q})$

History 2, which involves a discontinuous change at time $t=0$ from a 
total perimeter of 4 to a total perimeter of 6, is thus discontinuous
under $\mu$.

Over the space ${\mathcal D}$, $\Topo_{\mu}$ supports continuous
morphing; this is equivalent to saying that the perimeter is a continuous
function in the Hausdorff metric topology $\Topo_{H}$. Over the larger space
$\Reg$, 
$\Topo_{\mu}$ does not support continuous
morphing, as discussed above in section~\ref{secWellBehaved}.

In $\Euc^{n}$ for $n>2$, one might have pieces of any dimensionality $k<n$ 
peel off from the side:

{\bf History 3.k} ($k=0 \ldots n-1$):
In $\Euc^{n}$ 
let $\phi_{0} = (0,1)^{n}$ and let
$\phi(t) = (0,1)^{n} \cup(0,1)^{k} \times (1+t,1+2t)^{n-k}$.

The metric $H^{d}$ takes these all to be continuous. The metric
$M$ takes them all to be discontinuous. If one defines a metric 
$\mu({\bf P,Q})$ as the sum of $H_{d}({\bf P,Q})$ plus the absolute value 
of the difference of the kth order
quermassintegrals, then history 3.k will be continuous
for all $k<j$ and discontinuous for all $k\geq j$.

In ${\mathcal D}$,
histories such as 3.k for $k > 0$ can only be constructed 
starting if part of $\phi_{0}$ is a $k$-dimensional flat surface. If 
$\phi_{0}$ is strongly convex, then only the analogue of history 3.0 can
be constructed. Equivalently, over the space of regions whose closure is
strictly convex, the metrics defined above all define the same topology
for all values of $k$.

\subsection{The Hausdorff metric in ${\mathcal D}$}

\begin{figure}
\begin{center}
\includegraphics[width=4in]{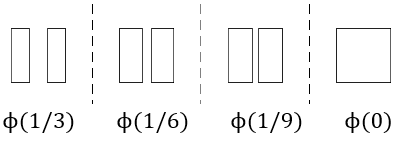}
\end{center}
\caption{History 4}
\label{figHistory4}
\end{figure}

The Hausdorff distance $H({\bf P,Q})$ is always greater than or equal to
the dual-Hausdorff distance; hence the topology it generates is coarser.
Indeed over the space $\mathcal D$
it is strictly coarser, as history 4 illustrates 
(figure~\ref{figHistory4})

{\bf History 4:} \\
$\phi(0) = (0,2) \times (0,2)$. \\
$\phi(t) = (0,1-t) \times (0,2) \cup (1+t,2) \times (0,2).$

For $t>0$, $H(\phi(t),\phi(0)) = t$; every point of $\phi(t)$ is in $\phi(0)$
and every point in $\phi(0)$ is within $t$ of $\phi(t)$. On the other hand
for all $t$ $H^{d}(\phi(t),\phi(0)) = 1$; the point $\la 1,1 \ra$ is
in $\phi(t)^{c}$ but is distance 1 from any point in $\phi(0)^{c}$.
Thus History 4 is continuous at time $t=0$ under the Hausdorff distance
but discontinuous over the dual-Hausdorff distance.

\subsection{The symmetric-difference metric in ${\mathcal D}$}

\begin{lemma}
\label{lemDilationHausdorff}
Let {\bf P} and {\bf Q} be regions such that $H^{1}({\bf Q,P}) \leq \delta$. 
Let ${\bf W}(\delta)$ be the dilation of {\bf P} by $\delta$.
Then ${\bf Q} \subset {\bf W}(\delta)$. 
\end{lemma}

{\bf Proof:} Immediate from the definitions.

\begin{lemma}
\label{lemConvexErosion}
Let {\bf P} and {\bf Q} be convex regions. Let $\delta > H({\bf P,Q})$
Then $\mbox{erode}({\bf P},\delta) \subset {\bf Q}$.
\end{lemma}

{\bf Proof:} of the contrapositive. Suppose that point 
${\bf x} \in \mbox{erode}({\bf P},\delta)$ and that
${\bf x} \not \in {\bf Q}$. Since $\bf Q$ is convex, there is a plane
$\bf Z$ through {\bf x} such that {\bf Q} lies on one side of {\bf Z}. Let 
$\bf H$ be the hemisphere of $\bar{B}({\bf x},\delta)$ on the far side
of {\bf Z} from {\bf Q} and let {\bf c} be the apex of {\bf H}. Then
${\bf c} \in \bar{{\bf P}}$ and $d({\bf c,Q}) \geq \delta$, so $H({\bf P,Q}) \geq
\delta$.

\begin{corollary}
\label{corSymDiffInShell}
If $\bf P$ and $\bf Q$ are convex then the symmetric difference of {\bf P}
and {\bf Q} is a subset of
the union of the inner and outer shells of {\bf P} by
the Hausdorff distance. \\
${\bf S}({\bf P,Q}) \subset 
{\bf O}({\bf P},H({\bf P,Q})) \cup 
{\bf I}({\bf P},H({\bf P,Q}))$
\end{corollary}

{\bf Proof:} Immediate from lemmas~\ref{lemDilationHausdorff}
and \ref{lemConvexErosion}.

\begin{lemma}
\label{lemShell}
Let $\bf P$ be any bounded open
region. Then for any $\epsilon > 0$, 
there exists $\delta > 0$ such that 
$v({\bf O}({\bf P},\delta)) < \epsilon$ and 
$v({\bf I}({\bf P},\delta)) < \epsilon$.
\end{lemma}

{\bf Proof:} Easily shown from the definition of measure as a limit.

\begin{lemma}
\label{lemSmallerShell}
Let {\bf P} be a convex region, let $\epsilon > 0$.
and let ${\bf Q}$ be a convex region such that
$\mbox{dilate}({\bf Q}, \epsilon) \subset {\bf P}$.
Then $v({\bf O}({\bf Q},\epsilon)) \leq v({\bf O}({\bf P},\epsilon))$.
\end{lemma}

{\bf Proof:} 
Let ${\bf X} = \mbox{dilate}({\bf Q},\epsilon)$. 
Let ${\bf Z} \subset {\bf X}$ be a convex polytope such that 
$v({\bf X} \setminus {\bf Z}) < \alpha$.

Let ${\bf Y}_{1} \ldots {\bf Y}_{m}$ be the faces of {\bf Z}.
For $i=1 \ldots m$: 
let ${\bf C}_{i}$ be the prism where one face is ${\bf Y}_{i}$, the
axis has length $\epsilon$, is orthogonal to ${\bf Y}_{i}$ and extends
inward into ${\bf Z}$.

I claim that 
$\bigcup_{i=1}^{m} {\bf C}_{i} \supset {\bf Z} \cap {\bf O}({\bf Q},\epsilon)$.
Proof: Let $\bf z$ be a point in ${\bf Z} \cap {\bf O}({\bf Q},\epsilon)$.
Let {\bf a} be the closest point to {\bf z} on $@{\bf X}$. 
Let {\bf b} be the intersection of the line {\bf az} with $@{\bf Z}$.
Let {\bf c} be the closest point to {\bf z} on $@{\bf Z}$. 
Let ${\bf Y}_{i}$ be the face of {\bf z} containing {\bf c}.
Then 
$\epsilon \geq d({\bf z,a}) \geq d({\bf z,b}) \geq d({\bf z,c})$. Moreover
the line {\bf zc} is orthogonal to ${\bf Y}_{i}$, so ${\bf z} \in {\bf C}_{i}$.

Therefore $v({\bf O}({\bf Q},\epsilon)) \leq 
v(\bigcup_{i=1}^{m} {\bf C}_{i}) + v({\bf X} \setminus {\bf Z}) \leq
v(\bigcup_{i=1}^{m} {\bf C}_{i}) + \alpha \leq 
\sum_{i=1}^{m} v({\bf C}_{i})$. 

Now extend each prism ${\bf C}_{i}$ outward from @{\bf Z}. Let ${\bf D}_{i}$
be the intersection of each such extended prism with 
${\bf O}({\bf P},\epsilon)$. Since {\bf Z} is convex, no two of these intersect.
Moreover, each ${\bf D}_{i}$ contains a right prism with cross section ${\bf Y}_{i}$
and with length at least $\epsilon$, so $v({\bf D}_{i}) \geq v({\bf C}_{i})$. \\
So $v({\bf O}({\bf P},\epsilon)) \geq
\sum_{i=1}^{m} v({\bf D}_{i}) \geq 
\sum_{i=1}^{m} v({\bf C}_{i}) \geq 
v({\bf O}({\bf Q},\epsilon)) - \alpha$. Since $\alpha$ can be made arbitrarily
small, we have
$v({\bf O}({\bf P},\epsilon)) \geq
v({\bf O}({\bf Q},\epsilon))$.

$\QED$

\begin{corollary}
\label{corUniformShell}
Let ${\bf P}$ be a convex region and let $\epsilon >0$. Then there exists
$\delta>0$ such that, for any convex region ${\bf Q} \subset {\bf P}$,
$v({\bf O}({\bf Q},\delta)) < \epsilon$.
\end{corollary}

{\bf Proof:} 
Choose $\delta_{1} > 0$. 
Let ${\bf W} = {\bf O}({\bf P},\delta_{1})$.
Using lemma~\ref{lemShell}, choose $\delta_{2}$ so that 
$v({\bf O}({\bf W},\delta_{2})) < \epsilon$. 
Let $\delta= \min(\delta_{1},\delta_{2})$.
Then since $\mbox{dilate}({\bf Q},\delta) \subset {\bf W}$, by
lemma~\ref{lemSmallerShell}, $v({\bf O}({\bf Q},\delta)) < \epsilon$.

\begin{theorem}
\label{thmVolumeCoarserThanHausdorff}
$\Topo_{V}$ is strictly coarser than $\Topo_{H}$ over $\mathcal D$.
\end{theorem}

{\bf Proof:} We first prove that $\Topo_{H}$ is at least as fine as
$\Topo_{V}$ over $\mathcal D$. We need to show that, for any region 
${\bf P} \in \mathcal D$  and $\epsilon > 0$ there exists $\delta > 0$ such
that, if ${\bf Q} \in \mathcal D$ and $H({\bf Q,P}) < \delta$ then 
$V({\bf Q,P}) < \epsilon$.

Choose {\bf P} and $\epsilon > 0$. There are two cases:

Case 1: {\bf P} is convex. 
Using lemma~\ref{lemShell} choose $\delta_{1}$ such that 
$v({\bf O}({\bf P},\delta_{1})) < \epsilon/4$ and
$v({\bf I}({\bf P},\delta_{1})) < \epsilon/4$.
Then, by corollary~\ref{corSymDiffInShell}
for every convex $\bf Q$, if $H({\bf P,Q}) < \delta_{1}$,  
$v({\bf S(P,Q)}) < \epsilon/2$.

Let ${\bf W} =\mbox{dilate}({\bf P},\delta_{1})$.
Using corollary~\ref{corUniformShell} choose $\delta_{2}$ such that,
for every convex subset {\bf X} of {\bf W}, 
$v({\bf O}({\bf X},\delta_{2}) < \epsilon/4$. 
Let $\delta =\min(\delta_{1},\delta_{2})$.

Suppose that ${\bf Q} \in {\mathcal D}^{2}$ such that $H({\bf Q,P}) < \delta$.
Let ${\bf Q} = {\bf C} \cup {\bf D}$ where ${\bf C}$ and ${\bf D}$ are convex.
Since $H^{1}({\bf Q,P}) < \delta$ it follows that ${\bf Q} \subset {\bf W}$.
Hence $v({\bf Q} \setminus {\bf P}) \leq
v{(\bf W} \setminus {\bf P}) \leq \epsilon/2$. 

Since $H^{1}({\bf P,Q}) < \delta$ it follows that 
${\bf P} \subset \mbox{dilate}({\bf Q},\delta) = 
\mbox{dilate}({\bf C},\delta) \cup \mbox{dilate}({\bf D},\delta)$. \\
Hence ${\bf P} \setminus {\bf Q} \subset 
(\mbox{dilate}({\bf C},\delta) \cup \mbox{dilate}({\bf D},\delta)) 
\setminus {\bf Q} \subset {\bf O}({\bf C},\delta) \cup {\bf O}({\bf D},\delta)$.
\\
But $\mbox{dilate}({\bf C},\delta)$ and $\mbox{dilate}({\bf D},\delta)$ are 
both convex subsets of {\bf W}, so
$v({\bf O}({\bf C},\delta) \leq \epsilon/4$ and
$v({\bf O}({\bf D},\delta) \leq \epsilon/4$.
So $v({\bf P} \setminus {\bf Q}) < \epsilon/2$ and 
$v({\bf S(P,Q)}) < \epsilon$.

Case 2:  ${\bf P} \in {\mathcal D}^{2}$. By 
lemma~\ref{lemNoConvexClosetoNonConvex} there exists $\delta_{1} > 0$
such that there are no convex regions ${\bf Q}$ with $H({\bf P,Q}) <
\delta_{1}$.

Let ${\bf P} = {\bf C} \cup {\bf D}$
where ${\bf C}$ and {\bf D} are convex. By lemma~\ref{lemMatchPieces}
there exists $\delta_{2}> 0$, such that, for any 
${\bf Q} \in {\mathcal D}^{2}$, if $H({\bf P,Q}) < \delta_{2}$ then,
{\bf Q} can be divided into convex components {\bf E} and {\bf F}
such that $H({\bf C,E}) < \delta_{2}$ and $H({\bf D,F}) < \delta_{2}$.
Clearly 
${\bf S}({\bf P,Q}) =  {\bf S}({\bf C,E}) \cup {\bf S}({\bf D,F})$.
Using theorem~\ref{thmConvex} one can choose $\delta_{3}$ 
such that, if $H({\bf C,E}) < \delta_{3}$ then
$v({\bf S}({\bf C,E})) < \epsilon/2$ and
$v({\bf S}({\bf D,F})) < \epsilon/2$. Thus if 
$H({\bf P,Q}) < \min(\delta_{1},\delta_{3})$ then 
$V({\bf P,Q}) < H({\bf P,Q})$.

To show that $\Topo_{V}$ is strictly coarser than $\Topo_{H}$, note that 
histories 5.1 and 5.2 below
are continuous in $\Topo_{V}$ but not in $\Topo_{H}$. In
history 5.1 for $t > 0$, $V(\phi(t),\phi(0)) = t^{2}$ while 
$H(\phi(t),\phi(0)) = 1+t$.

$\QED$

{\bf History 5.1:} \\
$\phi(0) = (0,1) \times (0,1).$ \\
$\phi(t) = (0,1) \times (0,1) \cup (2,2+t) \times (0,t)$ for $t > 0$.

{\bf History 5.2:} \\
$\phi(0) = (0,1) \times (0,1).$ \\
$\phi(t) = (0,1) \times (0,1) \cup (2,2+t) \times (0,1)$ for $t > 0$.

Analogous with histories 3.k, in $\Euc^{n}$, one can define $n$ 
qualitatively different histories, depending on the dimensionality of
the new piece.

{\bf History 6.k} ($k=0 \ldots n-1$)
In $\Euc^{n}$ 
let $\phi_{0} = (0,1)^{n}$ and let
$\phi(t) = (0,1)^{n} \cup(0,1)^{k} \times (2,2+2t)^{n-k}$.

As with histories 3.k, if one defines a metric 
$\mu({\bf P,Q})$ as the sum of $V({\bf P,Q})$ plus the absolute value 
of the difference of the $k$th-order
quermassintegrals, then history 6.k will be continuous
for all $k<j$ and discontinuous for all $k\geq j$. Unlike histories 3.k, these
multiple types of histories are possible even if $\phi_{0}$ is strictly
convex. (Define $\phi(t)$ as $\phi(0)$ union an ellipsoid with $k$ axes
of length 1 and $n-k$ axes of length $t$.)

\subsection{Wasserstein metrics in ${\mathcal D}$}

\begin{figure}
\begin{center}
\includegraphics[height=3in]{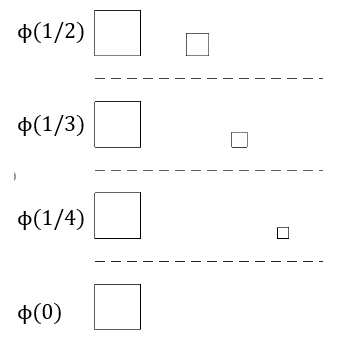}
\end{center}
\caption{History 7.$\psi$, where $\psi(t)=t^{2}$}
\label{figHistory7}
\end{figure}

To compare the topologies generated by the Wasserstein distances, we 
consider the following infinite collection of histories:

{\bf History 7}.$\psi$ (figure~\ref{figHistory7}). 
Let $\psi : \Real \mapsto \Real$ be a 
continuous function such that
$\psi(0)=0$ and $\lim_{x \goesto \infty} \psi(x) = \infty$. \\
Define the history $\phi^{\psi}: \Real \mapsto \Euc^{n}$ as: \\
$\phi^{\psi}(0)=(0,1)^{n}$. \\
$\phi^{\psi}(t) = \phi(0) \cup [(0,t)^{n-1} \times 
(\psi^{-1}(t^{-n}), \psi^{-1}(t^{-n})+t)]$.

The idea is that at time $t>0$, the unit box $(0,1)^{n}$ is joined by
another box of size $t^{n}$, 
growing from zero size, and heading inward 
from infinitely far away.
The trade-off between the size of the box and its distance is governed by the
function $\psi$ (the specific time dependence doesn't matter.)

\begin{lemma} 
\label{lemWasserstein1}
Let $\beta$ be a Mulholland functions. Let $\alpha(x)$ be a continuous
function such that
$\alpha(0)=0$ and $\lim_{x \goesto \infty} \alpha(x) = \infty$. \\
Let $\phi^{\alpha}(t)$ be as in History 7.$\alpha$. Then
\[ \lim_{t \goesto 0^{+}} W^{\beta}(\phi^{\alpha}(t),\phi^{\alpha}(0)) =
\left\{ \begin{array}{ll}
0 & \mbox{if } \lim_{x \goesto \infty} \beta(x)/\alpha(x) = 0 \\
\infty & \mbox{if } \lim_{x \goesto \infty} \beta(x)/\alpha(x) = \infty  
\end{array} \right. \]
\end{lemma}

{\bf Proof} (somewhat informal): 
The value of the integral in the definition of the Wasserstein
distance $W^{\beta}$ is dominated by the cost of moving the quantity $t^{n}$
of material a distance $d(t)=\alpha^{-1}(t^{-n}))$. By definition of the 
Wasserstein
distance, that cost 
$c(t) \approx \beta(d(t)) \cdot t^{-n} \approx 
\beta(\alpha^{-1}(t^{-n})) \cdot t^{n}$.
The Wasserstein distance is $W^{\beta}(\phi(0),\phi(t)) \approx 
\beta^{-1}(c(t)).$ So as $t \goesto \infty$, 
if $\beta(t) \ll \alpha(t)$,
then, as $t \goesto 0^{+}$,  $\beta(\alpha^{-1}(t^{-n})) \ll t^{-n}$ so
$c(t)$ and $W^{\beta}(t)$ go to 0;
if $\beta(t) \gg \alpha(t)$,
then, as $t \goesto 0^{+}$,  $\beta(\alpha^{-1}(t^{-n})) \gg t^{-n}$ so
$c(t)$ and $W^{\beta}(t)$ go to $\infty$.

\begin{lemma}
\label{lemWasserstein2}
Let $\alpha$, $\beta$  be two Mulholland functions. If
$\alpha(x) \ll \beta(x)$ as $x \goesto \infty$ then, 
over $\mathcal D$,
topology $\Topo_{{W}^{\alpha}}$ is not finer than the 
topology $\Topo_{{W}^{\beta}}$.
\end{lemma}

{\bf Proof:} Let $\zeta(x) = \sqrt{\alpha(x)\beta(x)}$
By lemma~\ref{lemWasserstein1} $\phi^\zeta(t)$ is continuous
relative to $\Topo_{W^{\alpha}}$ but discontinuous with respect to 
$\Topo_{W^{\beta}}$.

\begin{lemma}
\label{lemWasserstein3}
Let $\alpha$, $\beta$ be two Mulholland functions. If
$\alpha(x) \ll \beta(x)$ as $x \goesto \infty$ then, over $\Reg$,
topology $\Topo_{{W}^{\beta}}$ is at least as fine as
topology $\Topo_{{W}^{\alpha}}$.
\end{lemma}

{\bf Proof:} 
The intuition of the proof is this: Suppose that ${\bf Q}_{i}$ is close to 
$\bf P$ in the measure $W^{\beta}$. Let $\gamma$ be mapping of $\bf P$ to
${\bf Q}_{i}$ such that $C^{\beta}(\gamma,{\bf P},{\bf Q}_{i})$ is close
to $W^{\beta}({\bf P},{\bf Q}_{i})$.
Divide {\bf P} into two parts: the points that $\gamma$ is moving only a short
distance, and the points that it is moving a long distance. If you consider
now the integral using $\alpha$: the first part is moving only a small distance
so it makes a small contribution to the integral in $W^{\alpha}$. 
Over the second part, the integral using $\alpha$
can't be very much larger than the integral using $\beta$, so it is also
makes a small contribution to $W^{\alpha}$

Formally: We need to show that,  for any region {\bf P} and sequence
${\bf Q}_{1}, {\bf Q}_{2} \ldots$ if 
$W^{\beta}({\bf Q}_{i},{\bf P})$
converges to 0, then
$W^{\alpha}({\bf Q}_{i},{\bf P})$ also converges to 0.
Since $\alpha(x)$ and $\beta(x)$ go to 0 as $x$ goes to 0, in view
of the definition of $W^{\phi}$,
it clearly suffices to show that, for any $\epsilon > 0$ there exists
$\delta > 0$ such that, for any ${\bf Q}_{i}$ and uniform function $\gamma$
from  ${\bf P}$ to  ${\bf Q}_{i}$, 
if $I^{\beta}(\gamma,{\bf P}) < \delta$ then
$I^{\alpha}(\gamma,{\bf P}) < \epsilon$ where $I^{\psi}$ is the integral
defined earlier:

Choose $\epsilon > 0$. Let $M = \sup_{x \in [\alpha^{-1}(\epsilon/2),\infty)} 
\alpha(x)/\beta(x)$.  Since 
$\alpha(x) \ll \beta(x)$ as $x \goesto \infty$, this supremum exists
and is finite. Let $\delta = \epsilon/2M$. Assume that 
$I^{\beta}(\gamma,{\bf P}) < \delta$.
Partition {\bf P} into two subsets (either may be empty):
\\
${\bf P}_{1} = \{ {\bf x} \: | \: d({\bf x},\gamma({\bf x})) < 
\alpha^{-1}(\epsilon/2) \}$. \\
${\bf P}_{2} = \{ {\bf x} \: | \: d({\bf x},\gamma({\bf x}))  \geq
\alpha^{-1}(\epsilon/2) \}$. \\

Clearly 
\[ I^{\alpha}(\gamma,{\bf P}) =
\frac{1}{v({\bf P})} \int_{{\bf x} \in {\bf P}} 
\alpha(d({\bf x},\gamma({\bf x}))) \:
\mbox{d}{\bf x}  = 
 \frac{1}{v({\bf P})} \int_{{\bf x} \in {\bf P}_{1}} 
\alpha(d({\bf x},\gamma({\bf x}))) \:
\mbox{d}{\bf x}  +  
\frac{1}{v({\bf P})} \int_{{\bf x} \in {\bf P}_{2}} 
\alpha(d({\bf x},\gamma({\bf x}))) \:
\mbox{d}{\bf x}   \]

But for ${\bf x} \in {\bf P}_{1}$, 
$\alpha(d({\bf x},\gamma({\bf x}))) \leq \epsilon/2$, so
\[ \frac{1}{v({\bf P})} \int_{{\bf x} \in {\bf P}_{1}} 
\alpha(d({\bf x},\gamma({\bf x}))) \:
\mbox{d}{\bf x}  < \frac{v({\bf P}_{1})}{v({\bf P})}(\epsilon/2) \leq \epsilon/2 \]

And for ${\bf x} \in {\bf P}_{1}$, 
$\alpha(d({\bf x}, \gamma({\bf x}))) \leq 
M \beta(d({\bf x}, \gamma({\bf x})))]$ so
\[ 
\frac{1}{v({\bf P})} \int_{{\bf x} \in {\bf P}_{1}} 
\alpha(d({\bf x},\gamma({\bf x}))) \:
\mbox{d}{\bf x}  < 
\frac{1}{v({\bf P})} \int_{{\bf x} \in {\bf P}_{1}} 
M \beta(d({\bf x},\gamma({\bf x}))) \:
\mbox{d}{v} ) < 
M I^{\beta}(\gamma,{\bf P}) < \epsilon/2 \] 
 
$\QED$

\begin{theorem}
Let $\alpha$, $\beta$ be two Mulholland functions. If
$\alpha(x) \ll \beta(x)$ as $x \goesto \infty$ then, over $\Reg$ and over
$\mathcal D$,
topology $\Topo_{{W}^{\beta}}$ is strictly finer than
topology $\Topo_{{W}^{\alpha}}$.
\end{theorem}

{\bf Proof:} Immediate from lemmas~\ref{lemWasserstein2} and 
\ref{lemWasserstein3}. $\QED$

With a slight modification of the proof of \ref{lemWasserstein3} we can show 
that, if you consider a bounded subset of $\Reg$, then
any two Wasserstein distances give the identical topology. In other words
if you want to construct an example like History.7.$\psi$ that is continuous
relative to one Wasserstein distance and discontinuous relative to another,
then you have to use a similar construction of using, as $t \goesto 0^{+}$
smaller and smaller regions further and further out.

\begin{theorem}
\label{thmWassersteinBounded}
Let $\bf U$ be a bounded region in $\Euc^{n}$. Let $\mathcal V$ be
any collection of sub-regions of $\bf U$. 
Let $\alpha$, $\beta$ be two Mulholland functions. Then over 
$\mathcal U$, $\Topo_{W^{\alpha}} = \Topo_{W^{\beta}}$.
\end{theorem}

{\bf Sketch of proof:} Suppose that ${\bf Q}_{1}, {\bf Q}_{2} \ldots$
converges to {\bf P}, where these are all subsets of $\bf U$.
Suppose that this converges in $\beta$.
As in the proof of lemma~\ref{lemWasserstein3}, divide {\bf P} into
two parts; ${\bf P}_{1}$, the points that are being moved a short distance, 
and ${\bf P}_{2}$ the points that are being moved a long distance. The
integral over ${\bf P}_{1}$ is necessarily small using any Mulholland function.
Since the integral over ${\bf P}_{2}$ is small using $\beta$, and since 
the distance that points are being moved is not small, the volume of 
${\bf P}_{1}$ itself must be small. But the distance that they are being
moved cannot be more than diameter({\bf U}). Therefore the integrand is not
greater than $\alpha$(diameter({\bf U})), and since this is being taken
over a small volume, the result is also small.

As with histories 3, 5, and 6, one can add another parameter $k$, which
is the dimensionality of the new piece that appears.

{\bf History 7}.$\psi.k$: Let  $\psi : \Real \mapsto \Real$ be a 
continuous function such that
$\psi(0)=)$ and $\lim_{x \goesto \infty} \psi(x) = \infty$. Let $k$ 
be an integer between 0 and $n-1$\\
Define the history $\phi^{\psi}: \Real \mapsto \Euc^{n}$ as: \\
$\phi^{\psi}(0)=(0,1)^{n}$. \\
$\phi^{\psi}(t) = \phi(0) \cup [(0,1)^{k} \times (0,t)^{n-1} \times 
(\psi^{-1}(t^{k-n}), \psi^{-1}(t^{k-n})+t)]$.

It is easily seen that, 
if one considers metrics which are the sum of a Wasserstein function plus
the absolute value of the difference of the $k$th-order quermassintegrals,
then, for any two functions $\alpha$ and $\beta$ and any two values $k,m$
between 0 and $n-1$, if either $\alpha$ and $\beta$ have different growth
rates or $k \neq m$, then one can construct a history $\phi$ of this form
which is continuous with respect to one metric and discontinuous with 
respect to the other. Thus any two such metrics generate different topologies.
The distinction between different values of $k$ can be achieved even
if the space of regions is limited to subsets of a bounded region.

\begin{lemma}
\label{lemWassersteinHausdorffConvex2A}
Let $\alpha$ be a Mulholland function. Then over $\mathcal D$, the
corresponding Wasserstein metric topology $\Topo_{W^{\alpha}}$ 
generates a topology that
is not finer than the Hausdorff metric topology $\Topo_{H}$.
\end{lemma}

{\bf Proof:} Consider history 5.1 above: \\
$\phi(0) = (0,1) \times (0,1).$ \\
$\phi(t) = (0,1) \times (0,1) \cup (2,2+t) \times (0,t)$ for $t > 0$  \\

It is easily shown that $H(\phi(0),\phi(t)) = 1+t$ but for any $\alpha$,
$W^{\alpha}(\phi(0),\phi(t)) \approx \alpha^{-1}(t^{2})$. Thus, $\phi$ is continuous
at $t=0$ in the Wasserstein topology but discontinuous in the 
Hausdorff-metric topology.

\begin{lemma}
\label{lemHausdorffFinerWasserstein}
Let $\alpha$ be a Mulholland function. Then over $\mathcal D$, the
corresponding Wasserstein metric topology $\Topo_{W^{\alpha}}$ 
generates a topology that
is coarser than the Hausdorff metric topology $\Topo_{H}$.
\end{lemma}

{\bf Proof:} Choose region ${\bf P} \in {\mathcal D}$ 
and $\epsilon > 0$. 
Let $p=\mbox{diameter}({\bf P})$. Let $a=\alpha(\epsilon)v({\bf P})/2p$.
Using theorem~\ref{thmVolumeCoarserThanHausdorff}, choose $\delta_{1}$
such that, for all ${\bf Q} \in {\mathcal D}$, if 
$H({\bf P,Q}) < \delta_{1}$ then $V({\bf P,Q}) < a$.
Let $\delta= \min(\delta_{1},p/2)$. Then by 
lemma~\ref{lemTopoDHFinerWasserstein2},
if ${\bf Q} \in {\mathcal D}$ and $H({\bf P,Q}) < \delta$, then
$W^{\psi}({\bf P,Q}) < \epsilon$.

\begin{theorem}
\label{thmHausdorffFinerWasserstein}
Over the space $\mathcal D$, the Hausdorff metric
topology is strictly finer than any Wasserstein metric topology.
\end{theorem}

{\bf Proof:} This is the combination of 
lemmas~\ref{lemWassersteinHausdorffConvex2A} and`
\ref{lemHausdorffFinerWasserstein}.

\begin{lemma}
\label{lemSymDifNotFiner}
Over $\mathcal D$, the symmetric difference topology $\Topo_{V}$
is not finer than any Wasserstein metric topology
$\Topo_{W^{\alpha}}$.
\end{lemma}

{\bf Proof:} Let $\psi({\bf x}) = \alpha^{2}({\bf x})$. Then History.7.$\psi$
is continuous relative to $\Topo_{V}$ 
but not relative to $\Topo_{W^{\alpha}}$ 
by lemma~\ref{lemWasserstein1}.

\begin{lemma} 
\label{lemDiffOfConvexA}
Let {\bf P} be a convex region and $\epsilon > 0$. Then there
exists $\delta > 0$ such that, for any convex {\bf Q}, if 
$v({\bf P} \setminus {\bf Q}) > \epsilon$ then there exists a point {\bf p}
such that ${\bf B}({\bf p},\delta) \subset {\bf P} \setminus {\bf Q}$.
\end{lemma}

{\bf Proof:} Using lemma~\ref{lemShell}, choose $\delta_{1}$ such that 
$v({\bf I}({\bf P},\delta_{1})) < \epsilon$.
Let ${\bf R} = \mbox{erode}({\bf P},\delta_{1}).$
Let ${\bf Q}$ be a convex region such that
$v({\bf P} \setminus {\bf Q}) > \epsilon$. Clearly $\bf R$ is not a subset
of ${\bf Q}$ since $v({\bf P} \setminus {\bf R}) < \epsilon$. Let ${\bf r}$ be a 
point in {\bf R} but not in {\bf Q}. Since ${\bf r} \in {\bf R}$ it follows
that ${\bf B}({\bf r},\delta_{1}) \subset {\bf P}$; since {\bf Q} is convex,
there is
at least a hemisphere of ${\bf B}({\bf r},\delta_{1})$ that is not in {\bf Q}.
Therefore there is a ball of radius $\delta_{1}/2$ in 
${\bf P} \setminus {\bf Q}$.

\begin{lemma} 
\label{lemDiffOfConvexB}
Let {\bf P} be a convex region and $\epsilon > 0$. Then there
exists $\delta > 0$ such that, for any ${\bf Q} \in {\mathcal D}^{2}$, 
if $v({\bf P} \setminus {\bf Q}) > \epsilon$ then there exists a point {\bf p}
such that ${\bf B}({\bf p},\delta) \subset {\bf P} \setminus {\bf Q}$.
\end{lemma}

{\bf Proof:} Choose $\bf P$ and $\epsilon > 0$. 
Using corollary~\ref{corUniformShell} choose 
$\delta_{1} > 0$ such that, for all convex
${\bf X} \subset {\bf P}$, $v({\bf O}({\bf X},\delta_{1})) < \epsilon/2$.
Let $\delta = \delta_{1}/2$.

Let ${\bf Q}$ be any region in ${\mathcal D}^{2}$ 
such that $v({\bf P} \setminus {\bf Q}) > \epsilon$.
Let {\bf C} and {\bf D}
be the two components of {\bf Q}. Let 
${\bf C}' = {\bf C} \cap {\bf P}$ and
${\bf D}' = {\bf D} \cap {\bf P}$. If either of these is empty, then the
result follows from lemma~\ref{lemDiffOfConvexA}, so assume that neither is
empty. Let ${\bf X}$ be a hyperplane dividing ${\bf C}'$ from ${\bf D}'$. Then
{\bf X} divides ${\bf P}$ into two parts, $\bf E$ containing $\bf C$ and
{\bf F} containing {\bf D}. 

Clearly {\bf E} and {\bf F} are convex and
${\bf P} \setminus {\bf Q} = 
({\bf E} \setminus {\bf C}) \cup 
({\bf F} \setminus {\bf D})$. Therefore either
$v({\bf E} \setminus {\bf C}) >  \epsilon/2$ or
$v({\bf F} \setminus {\bf D}) >  \epsilon/2$.  Assume the former. By the same
argument as in lemma~\ref{lemDiffOfConvexA}, there exists a point 
${\bf r}$ such that 
${\bf B}({\bf r},\delta) \subset {\bf E} \setminus {\bf C}$.

\begin{lemma} 
\label{lemDiffOfConvexC}
Let {\bf P} be a region in $\mathcal D$ and $\epsilon > 0$. Then there
exists $\delta_{1}, \delta_{2} > 0, $ such that, for any region {\bf Q}, 
if $v({\bf Q} \setminus {\bf P}) > \epsilon$ then there is a subset
${\bf W} \subset {\bf Q}$ such that $d({\bf P,W}) > \delta_{1}$ and
$v({\bf W})/v({\bf Q}) > \delta_{2}$.
\end{lemma}

{\bf Proof:} Using lemma~\ref{lemShell}, choose $\delta_{1}$ such that 
$v({\bf O}({\bf P},\delta_{1})) < \epsilon/2$.
Let ${\bf R} = \mbox{dilate}({\bf P},\delta_{1}).$
Let ${\bf Q}$ be a region such that $v({\bf Q} \setminus {\bf P}) > \epsilon$. 
Let ${\bf W} = {\bf Q} \setminus {\bf R}$. Then
${\bf Q} \setminus {\bf P} \subset {\bf W} \cup ({\bf R} \setminus {\bf P})$
so $v({\bf W}) > \epsilon/2$. So the conclusion is satisfied with
$\delta_{2} = \epsilon/(\epsilon + v({\bf R}))$.

\begin{lemma} 
\label{lemDiffOfConvexD}
Let {\bf P} be a region in $\mathcal D$ and $\epsilon > 0$. Then there
exists $\delta_{1}, \delta_{2} > 0, $ such that, for any region 
${\bf Q} \in {\mathcal D}$, 
if $V({\bf P,Q}) > \epsilon$ then there is a subset
${\bf W} \subset {\bf Q}$ such that $d({\bf P,W}) > \delta_{1}$ and
$v({\bf W})/v({\bf Q}) > \delta_{2}$.
\end{lemma}

{\bf Proof:}
$V({\bf P,Q}) = v(({\bf P} \setminus {\bf Q}) \cup ({\bf Q} \setminus {\bf P}))$,
so if $V({\bf P,Q}) > \epsilon$ then either 
$v({\bf P} \setminus {\bf Q}) > \epsilon/2$ or
$v({\bf Q} \setminus {\bf P}) > \epsilon/2$.

Using lemmas~\ref{lemDiffOfConvexA} and \ref{lemDiffOfConvexB}, we can find
$\delta_{A}$ such that, for all ${\bf Q} \in \mathcal D$, if 
$v({\bf Q} \setminus {\bf P}) < \epsilon/2$ and
$v({\bf P} \setminus {\bf Q}) > \epsilon/2$, then there is a point {\bf r}
such that ${\bf B}({\bf r},\delta_{A}) \subset {\bf Q} \setminus {\bf P}$,
so in this case, we can choose ${\bf W}={\bf B}({\bf r},\delta_{A}/2)$.
Let $s= v({\bf B}({\bf r},\delta_{A}))$, the volume of the $n$-dimensional
sphere of radius $\delta_{A}$. 
Then $v({\bf W})/v({\bf Q}) \geq s/(v({\bf P})+\epsilon/2)$.

Using lemma~\ref{lemDiffOfConvexB} we can find 
$\delta_{B}, \delta_{C}$ such that, for all regions ${\bf Q}$, if
$v({\bf Q} \setminus {\bf P}) > \epsilon/2$ 
then there exists a subset ${\bf W} \subset {\bf Q} \setminus {\bf P}$ 
such that $d({\bf W,P}) > \delta_{B}$ and $v({\bf W})/v({\bf Q}) > \delta_{C}$.

So if we take $\delta_{1} = \min(\delta_{A}/2,\delta_{B})$ and
$\delta_{2} = \min(s/(v({\bf P})+\epsilon/2),\delta_{C})$,
the conclusion of the lemma is satisfied.

\begin{lemma}
\label{lemWassersteinFinerSymDiff}
Over $\mathcal D$, the symmetric difference topology $\Topo_{V}$
is coarser than any Wasserstein metric topology
$\Topo_{W^{\psi}}$.
\end{lemma}

{\bf Proof:} We need to show that, for any Mulholland function $\psi$,
for any ${\bf P} \in {\mathcal D}$ and $\epsilon > 0$ there exists $\delta >0$ 
such that, for any ${\bf Q} \in {\mathcal D}$, if 
$W^{\psi}({\bf P,Q}) < \delta$ then $V({\bf P,Q}) < \epsilon$.

Given $\psi, {\bf P}, \epsilon$ as above, by lemma~\ref{lemDiffOfConvexD}
there exist $\delta_{1}, \delta_{2}$ such that, for all ${\bf Q} \in 
\mathcal D$, if $V({\bf P,Q}) > \epsilon$ then there exists a region
${\bf W} \subset {\bf Q}$ 
such that $d({\bf W,P}) > \delta_{1}$ and $v({\bf W}) > \delta_{2}$.

Let $\gamma$ be any uniform mapping from $\bf Q$ to $\bf P$. Then
\[ I^{\psi}(\gamma) =
\int_{{\bf x} \in \bf Q} \psi(d({\bf x},\gamma({\bf x}))) \: d{\bf x} >
\int_{{\bf x} \in \bf W} \psi(d({\bf x},\gamma({\bf x}))) \: d{\bf x} >
\int_{{\bf x} \in \bf W} \psi(\delta_{1}) \: d{\bf x} >
\delta_{2}v({\bf Q}) \psi(\delta_{1}) \]

So $W^{\psi}({\bf P,Q}) = \inf_{\gamma} \psi^{-1}(1/v({\bf Q})) I(\gamma)
> \psi^{-1}(\delta_{2}\psi(\delta_{1}))$.

So the conclusion is satisfied with 
$\delta =\psi^{-1}(\delta_{2}\psi(\delta_{1}))$.
$\QED$

\begin{theorem}
\label{thmWassersteinFinerSymDiff}
Over the space ${\mathcal D}$,  any Wasserstein-metric topology is 
strictly finer than the symmetric-difference-metric topology.
\end{theorem}

{\bf Proof:} From lemmas \ref{lemSymDifNotFiner} and 
\ref{lemWassersteinFinerSymDiff}.

\section{Star-shaped regions}
\label{secStar}

Over the space $\mathcal S$ of star-shaped regions centered at the origin, 
the situation is very 
different. As we shall show, the Hausdorff metric, the Wasserstein metrics, 
and the symmetric difference metrics all yield topologies that are incomparable 
in terms of fineness.  

For simplicity, we will demonstrate our results in $\Euc^{2}$, but the 
generalizations to $\Euc^{n}, n>2$ are obvious. 
It will be convenient to define a generalized wedge function:

\begin{definition}
Let $\theta \in [0,2 \pi)$, $\beta \in (0,\pi/4), b \in (0,1),
l \in (0, \infty)$  The {\em wedge centered at $\alpha$ of width $\beta$ 
with base $b$ and length $l$,} denoted ${\bf G}(\alpha,\beta,b,l)$ is the
set of all points whose polar coordinate $\la r,\theta \ra$ 
satisfy $b <r < l$,
$\alpha-\beta/2 < \theta <  \alpha+\beta/2$.
\end{definition}

Note that $v({\bf G}(\alpha,\beta,b,l)) =
(l^{2}-b^{2})\beta.$

\begin{theorem}
Over $\mathcal S$, the symmetric-difference metric and the
Wasserstein metrics are not finer than the Hausdorff metric.
\end{theorem}

{\bf Proof:} Consider the following history $\phi(t)$: \\
{\bf History}.8 \\
$\phi(0) = {\bf B}(\vec{0},1)$. \\
$\phi(t) = {\bf B}(\vec{0},1) \cup {\bf G}(0,t,1,2)$. 

Then $H(\phi(t),\phi(0)) = 1$. 
$V(\phi(t),\phi(0)) = t$. It is easily to show, using lemma~\ref{lemProbDist},
that for any $\psi$, $\lim_{t \goesto 0^{+}} \: W^{\psi}(\phi(t),\phi(0))=0$.
Thus $\phi$ is continuous with respect to $V$ and to $W^{\psi}$ but not
with respect to $H$.

\begin{theorem}
\label{thmStarHausdorffSymDiff}
Over $\mathcal S$, the Hausdorff metric is not finer than 
the symmetric-difference metric and the Wasserstein metric.
\end{theorem}

\begin{figure}
\begin{center}
\includegraphics[width=4in]{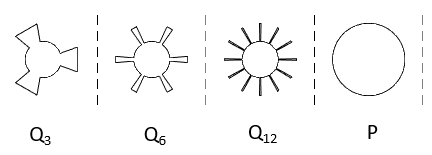}
\end{center}
\caption{Proof of theorem~\ref{thmStarHausdorffSymDiff}}
\label{figPorcupine1}
\end{figure}

{\bf Proof:} 
Let ${\bf P} = {\bf B}(\vec{0},2)$. For $k=1,2 \ldots$ 
let ${\bf Q}_{k} = {\bf B}(\vec{0},1) \cup 
\bigcup_{i=1}^{k} {\bf G}(2 \pi i/k, 2*\pi/k^{2},1,2)$ 
(Figure~\ref{figPorcupine1}).
That is, ${\bf Q}_{k}$ is the unit ball plus $k$ evenly spaced wedges of
width $1/k^{2}$ in the annulus between radius 1 and radius 2. As $k$ goes to 
infinity, the wedges get denser and denser within the ball of radius 2,
but the total area of the wedges is $6 \pi/k$. Thus $H({\bf P,Q}_{k}) \approx
1/2k - 1/2k^{2}$ but $V({\bf P,Q}) = 3\pi - 3/k$. Thus the sequence 
${\bf Q}_{k}$ converges to $\bf P$ with respect to the Hausdorff metric
but not with respect to the symmetric-difference metric.

To show that $W^{\psi}({\bf Q}_{k},{\bf P})$ does not converge to 0, note
that the fraction of the area of ${\bf Q}_{k}$ that is in the central ball
is $\pi/(\pi+3/k)$. Thus as $k \goesto \infty$, any uniform function
$\gamma$ from ${\bf Q}_{k}$ to {\bf P} must essentially spread the central
ball out over all of {\bf P}; the wedges become increasingly irrelevant.
So $\lim_{k \goesto \infty} 
W^{\psi}({\bf Q}_{k},{\bf P}) =
W^{\psi}({\bf B}(\vec{0},1),{\bf P})$.

\begin{theorem}
\label{thmStarHausdorffWasserstein}
Over $\mathcal S$, no Wasserstein metric is finer than the symmetric-difference
metric.
\end{theorem}

{\bf Proof:} We modify the example from the proof of 
theorem~\ref{thmStarHausdorffSymDiff} by making the central circle much 
smaller than the wedges.

\begin{figure}
\begin{center}
\includegraphics[width=4in]{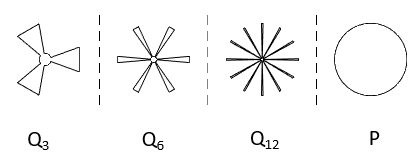}
\end{center}
\caption{Proof of theorem~\ref{thmStarHausdorffWasserstein}}
\label{figPorcupine2}
\end{figure}

Let ${\bf P} = {\bf B}(\vec{0},2)$. For $k=1,2 \ldots$ 
let ${\bf Q}_{k} = {\bf B}(\vec{0},1/k) \cup 
\bigcup_{i=1}^{k} {\bf G}(2 \pi i/k, 2\pi/k^{2},1/k,2)$.
(Figure~\ref{figPorcupine2}).

The combined area of the wedges approaches $4/k$, while the area of the central
circle is $\pi/k^{2}$. Define the mapping $\gamma$ from ${\bf Q}_{k}$
to {\bf P} so that, on the center circles $\gamma$ is the identity, and, 
on the edges, $\gamma$ spreads out the wedges uniformly in concentric circles
so that the entire
circle {\bf P} is covered.

\begin{quote}
For ${\bf x} \in {\bf B}(\vec{0},1/k)$, $\gamma({\bf x})={\bf x}$

For ${\bf x} \in {\bf G}(2 \pi i/k, 1/k^{2},1/k,2)$ if ${\bf X}$ has
polar coordinates
$\la r,\theta \ra$, then $\gamma({\bf x})$ has polar coordinates $\la r,
2 \pi i k + k \pi (\theta - 2 \pi i/k) \ra$.
\end{quote}

Let $\Gamma({\bf x})$ be the distribution generated by $\gamma$.
Almost all the mass in ${\bf Q}_{k}$ is in the wedges; in $\Gamma$ this mass
is distributed
evenly over the annulus $1/k < r < 2$. The density of $\Gamma$ over 
the inner circle
${\bf B}(\vec{0},1/k)$ is much larger, but that circle is small, so the total
mass there is small.
Therefore
using lemma~\ref{lemProbDist}, the distribution $\Gamma$
is close in Wasserstein distance to $U_{\bf Q}$.
However, $\gamma$ moves each
point by a maximum distance $2/k$; hence $W^{\psi}(U_{P},\Gamma)$ is small.
So for every $\psi$, 
$W^{\psi}({\bf Q}_{k},{\bf P})$ converges to 0 as $k \goesto \infty$.
However, $V({\bf Q}_{k},{\bf P}) = 4\pi - (3/k+ \pi/k^{2})$.

To compare Wasserstein functions over $\mathcal S$, we define a history
analogous to {\bf History}.7.$\psi$.

\begin{figure}
\begin{center}
\includegraphics[height=3in]{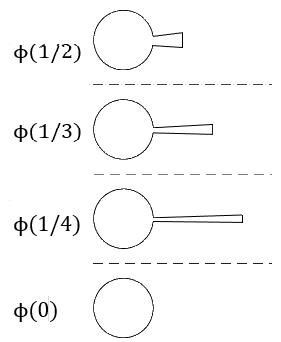}
\end{center}
\caption{History 9.$\psi$, with $\psi(t)=|t|$}
\label{figHistory9}
\end{figure}

{\bf History}.9.$\psi$. Let 
Let $\psi(x)$ be a continuous
function such that
$\alpha(0)=0$ and $\lim_{x \goesto \infty} \psi(x) = \infty$. \\
Let $\zeta$ be the inverse of $\psi$.
Define the history $\phi^{\psi}(t)$ as follows:

$\phi^{\psi}(0) = {\bf B}(\vec{0},1)$. \\
$\phi^{\psi}(t) = {\bf B}(\vec{0},1) \cup 
{\bf G}(0,t/\zeta^{2}(1/t),1,\zeta(1/t))$ (figure~\ref{figHistory9}).

\begin{lemma} 
\label{lemWassersteinStar1}
Let $\beta$ be a Mulholland functions. 
Let $\alpha(x)$ be a continuous
function such that
$\alpha(0)=0$ and $\lim_{x \goesto \infty} \alpha(x) = \infty$. \\
Let $\phi^{\alpha}(t)$ be as in History 9.$\alpha$. Then
\[ \lim_{t \goesto 0^{+}} W^{\beta}(\phi^{\alpha}(t),\phi^{\alpha}(0)) =
\left\{ \begin{array}{ll}
0 & \mbox{if } \lim_{x \goesto \infty} \beta(x)/\alpha(x) = 0 \\
\infty & \mbox{if } \lim_{x \goesto \infty} \beta(x)/\alpha(x) = \infty  
\end{array} \right. \]
\end{lemma}

{\bf Proof:} (Informal, analogous to the proof of lemma~\ref{lemWasserstein1}.)
A function $\gamma_{t}({\bf x})$ that transforms $\phi(0)$ into $\phi(t)$ 
involves, to order of magnitude, moving
a total of $t$ mass a distance of $\alpha^{-1}(1/t)$. Therefore the integral
$I(\gamma_{t})$  is roughly $t \cdot \beta(\alpha^{-1}(1/t))$. 
The Wasserstein distance is $W^{\beta}(\phi(0),\phi(t)) \approx 
\beta^{-1}()).$ So as $t \goesto \infty$, 
if $\beta(x) \ll \alpha(x)$,
then, as $t \goesto 0^{+}$,  $\beta(\alpha^{-1}(1/t)) \ll 1/t$ so
$I(\gamma_{t}$ and $W^{\beta}(t)$ go to 0;
if $\beta(t) \gg \alpha(t)$,
then, as $t \goesto 0^{+}$,  $\beta(\alpha^{-1}(t)) \gg t$ so
$I(\gamma_{t}$ and $W^{\beta}(t)$ go to $\infty$.

\begin{lemma}
\label{lemStarWasserstein}
Let $\beta$ be a Mulholland function and 
$\alpha(x) \ll \beta(x)$ as $x \goesto \infty$.  Then over $\mathcal S$,
$W^{\alpha}$ is not finer than $W^{\beta}$.
\end{lemma}

{\bf Proof:} Let $\zeta(x) = \sqrt{\alpha(x)\beta(x)}$
By lemma~\ref{lemWassersteinStar1} $\phi^\zeta(t)$ is continuous
relative to $\Topo_{W^{\alpha}}$ but discontinuous with respect to 
$\Topo_{W^{\beta}}$.

\begin{theorem}
\label{thmStarWasserstein}
Let $\beta$ be a Mulholland function and 
$\alpha(x) \ll \beta(x)$ as $x \goesto \infty$.  Then over $\mathcal S$,
$W^{\beta}$ is strictly finer than $W^{\alpha}$.
\end{theorem}

{\bf Proof:} Immediate from lemmas~\ref{lemWasserstein3} and 
\ref{lemStarWasserstein}.

\begin{theorem}
\label{thmStarWassersteinSymDiff}
Over $\mathcal S$, for any Mulholland function $\beta$, the 
symmetric-distance metric is not finer than the Wasserstein metric
$W^{\alpha}$.
\end{theorem}

{\bf Proof:} Using lemma~\ref{lemWassersteinStar1}, if 
$\psi = \sqrt{\alpha}$ then the function $\phi^{\psi}$ defined
in history.8.$\psi$ is continuous relative to the symmetric-difference
metric but not with respect to the metric $W^{\alpha}$.

\subsection*{Acknowledgements}
Thanks to Giorgio Stefani for helpful information; in particular, for drawing
my attention to the quermassintegral as a useful measure.

\end{document}